\documentclass[11pt]{amsart}
\usepackage{amsmath,amsthm, amscd, amssymb, amsfonts,mathtools}
\usepackage[all]{xy}
\usepackage{enumitem}

\usepackage{t1enc}
\usepackage[latin1]{inputenc}

\usepackage{fontsmpl}
\usepackage{srcltx}

\allowdisplaybreaks

\newcommand{\toba}{{\mathcal B}}
\newcommand{\fd}{finite-dimensional}

\usepackage{color}

\numberwithin{equation}{section}\theoremstyle{plain}

\newtheorem{theorem}{Theorem}[section]

\newtheorem{pro}[theorem]{Proposition}

\theoremstyle{definition}
\newtheorem{definition}[theorem]{Definition}

\theoremstyle{remark}
\newtheorem{rem}[theorem]{Remark}

\newcommand{\ydh}{{}^{H}_{H}\mathcal{YD}}

\def\G{\mathbb{G}}

\newcommand\ev{\operatorname{ev}}
\newcommand\id{\operatorname{id}}

\newcommand\ord{\operatorname{ord}}

\newcommand\End{\operatorname{End}}

\newcommand\GK{\operatorname{GK-dim}}

\def\q{\mathbf{q}}

\def\k{\Bbbk}
\def\ku{\Bbbk}

\def\ot{\otimes}

\def\N{\mathbb{N}}

\def\Cg{\mathfrak{C}}

\def\eps{\epsilon}

\def\bB{\mathfrak{B}}
\def\hit{\mathfrak{R}}
\def\Vhit{\mathbf{V}}

\def\ta{\texttt{a}}
\def\tb{\texttt{b}}

\newcommand{\J}{{\mathcal J}}

\newcommand\I{\mathbb I}

\def\lg{\langle}
\def\rg{\rangle}

\def\pf{\begin{proof}}
\def\epf{\end{proof}}

\begin{document}


 \title[Nichols algebras that are quantum planes]{Nichols algebras that are quantum planes}
\author[Andruskiewitsch and Giraldi]{Nicol\'as Andruskiewitsch and Jo\~ao Matheus Jury Giraldi}

\address{\noindent N. A. : FaMAF-CIEM (CONICET), Universidad Nacional de C\'ordoba,
Medina A\-llen\-de s/n, Ciudad Universitaria, 5000 C\' ordoba, Rep\'
ublica Argentina.}

\address{\noindent J. M. J. G. : Instituto de Matem\'atica, 
Universidade Federal do Rio Grande do Sul,
Rio Grande do Sul, Brazil.}

\email{andrus@famaf.unc.edu.ar,joaomjg@gmail.com}

\thanks{The work was partially supported by CONICET, Secyt (UNC), the MathAmSud project
GR2HOPF, CNPq (Brazil)}

\begin{abstract}
We compute all Nichols algebras of rigid vector spaces of dimension 2 that admit a non-trivial quadratic relation.
\end{abstract}
\maketitle
\section{Introduction}\subsection{}
Let $V$ be a vector space. We say that $c \in \End ( V \otimes V)$ satisfies the braid equation if
\begin{equation}\label{eqn:braid-bis}
(c\otimes \id)(\id \otimes c)(c\otimes \id) =  (\id \otimes c)(c\otimes \id)(\id \otimes c).
\end{equation}
If in addition $c$ is invertible, then we say that $(V,c)$ is a \emph{braided vector space} and that $c$ is a braiding. 
Assume that $\dim V < \infty$ and let $(v_i)$ be a basis of $ V $ and $(v^i)$ its dual basis. Then $c$ (or else $V$) is \emph{rigid} 
 if the  map $c^{\flat}: V^* \ot V \to V \ot V^*$ given by $ f\ot v \longmapsto \sum_{i}(\ev\ot \id \ot \id)(f\ot c(v\ot v_i)\ot v^i)$
 is invertible.

\smallbreak
Let $(V,c)$ be a rigid braided vector space and let $\toba(V)$ be its Nichols algebra, see \S \ref{subsec:nichols} for details.
An important problem is to determine the presentation and structure of $\toba(V)$, particularly
when it has finite dimension, or Gelfand-Kirillov dimension.
For example, let $n\in\mathbb N$ and $\I= \I_n = \{1,2,\dots,n\}$; let
 $\q = (q_{ij})_{i,j \in \I} \in (\ku^{\times})^{\I \times \I}$; let $V$ be a vector space with a basis $(v_{i})_{i \in \I}$; and let 
$c \in \End ( V \otimes V)$ be given by 
\begin{align*}
c(v_i\ot v_j) &= q_{ij} v_j\ot v_i, & i,j &\in \I.
\end{align*}
Then $(V,c)$ is a rigid braided vector space, called of \emph{diagonal type}.
The classification of the \fd{} Nichols algebras of diagonal type is known \cite{H-classif RS}. 
The following particular case was studied in \cite{AS1}. Assume that 
$q_{ij}q_{ji} =1$ for all $i \neq j \in \I$ and set $ N_i = \begin{cases} \ord q_{ii}, &\text{if } q_{ii} \neq 1; \\ \infty,  &\text{if } q_{ii} = 1.
\end{cases}$, $i\in \I$. Then
$\toba (V)$ is presented by generators $ x_{i}$,  $i \in \I$, with relations
\begin{align}\label{eq:qls1}
&& x_ix_j &= q_{ij} x_j x_i, & i &< j,
\\ \label{eq:qls2}
&& x_i^{N_i} &= 0, & \text{if }  N_i&< \infty.
\end{align}
In particular $\dim \toba(V) = \prod_{i\in \I} N_i$ is finite if and only if all $N_i$ are finite.  Also, $\GK \toba(V) = \vert\{i\in \I: N_i = \infty \} \vert$.

\smallbreak 
An algebra $A$ presented by generators $ x_{i}$,  $i \in \I$, with  relations \eqref{eq:qls1}, \eqref{eq:qls2}  is known as a quantum linear space,
or a \emph{quantum plane} when $n = 2$. 
It is well-known  that  for any quantum linear space $A$, we have
\begin{align}
\label{eq:qls3}
\{x_1^{a_1} \dots x_n^{a_n}: & \ 0 \le a_i < N_i,\, i\in \I \}
& \text{is a basis of } & A.
\end{align}

\subsection{}
In \cite[Propositions 4.8, 4.9]{GGi}, there were found braided vector spaces $(V,c)$ of dimension 2 not of diagonal type but such that $\toba(V)$
are quantum planes with $N_1= 2$ and $N_2 = 4$.
Notice that the coproduct of a Nichols algebra is determined by the braiding $c$, 
and vice versa the algebra and coalgebra structures determine the braiding \cite{Sch-ny};
thus these Nichols algebras are not isomorphic to 
Nichols algebras of diagonal type as a braided Hopf algebras. 
Consequently, the following question arises naturally: classify all Nichols algebras (of rigid braided vector spaces) that are isomorphic to
quantum linear spaces as algebras. In the present paper we solve this question for quantum planes; we recover 
the examples in \cite{GGi} in Remarks \ref{rem:GGi-uno} and \ref{rem:GGi-dos}.

\smallbreak
More generally, we consider braided vector spaces of
dimension 2 and compute the quadratic relations of their Nichols algebras, which is straightforward. 
Roughly speaking, there are four possible scenarios:
\begin{enumerate}[leftmargin=*,label=\rm{(\Alph*)}]
	\item There are no quadratic relations. We  plan to address (at least some of) these in the future.
	
\smallbreak	\item There is a quadratic relation close enough to \eqref{eq:qls1} that, together with \eqref{eq:qls2} for suitable $N_1$ and $N_2$, is a presentation of $\toba(V)$.
	Furthermore there is a PBW-basis \eqref{eq:qls3}.
	
\smallbreak	\item There are strange examples with quadratic relations and  more relations in higher degree  with  basis of a very specific type.
	
\smallbreak	\item There is a unique quadratic relation $x_1^2 = 0$ and an extra cubic relation that together present $\toba(V)$. The super Jordan plane \cite{AAH} fits here. 
\end{enumerate}

\begin{table}[ht]
	\caption{Nichols algebras of rank 2 }\label{tab:general}
	\begin{center}
		\begin{tabular}{| c | c | p{2.16cm}| c |p{2cm} |} 
			Name	& $R$-matrix & $\toba(V)$ & \small{Ref.}  & \small{Conditions}  	\\\hline
			\hline
			$\hit_{2, 1}$  & $\begin{pmatrix}
			k^2 & 0 & 0 & 0 \\ 
			0 & kp & k^2-pq & 0 \\ 
			0 & 0 & kq & 0 \\ 
			0 & 0 & 0 & k^2
			\end{pmatrix}$ 
			& \eqref{eqn:Nichols2,1} & \S \ref{subsec:hit21} & $k^2 =-1$ or $pq =1$  \\ \hline
			$\hit_{2, 2}$ & $\begin{pmatrix}
			k^2 & 0 & 0 & 0 \\ 
			0 & kp & k^2-pq & 0 \\ 
			0 & 0 & kq & 0 \\ 
			0 & 0 & 0 & -pq
			\end{pmatrix}$ & \eqref{eqn:Nichols2,2}, \eqref{eqn:Nichols2,2a} & \S \ref{subsec:hit22} & $k^2 =-1$ or $pq =1$ 
			\\ \hline
			$\hit_{2, 3}$ & $\begin{pmatrix}
			k & p & q & s \\ 
			0 & k & 0 & q \\ 
			0 & 0 & k & p \\ 
			0 & 0 & 0 & k
			\end{pmatrix}$
			& Table \ref{tab:h23} & \S \ref{subsec:hit23} & $ k^2=1 $ \\		
			\hline
			$\hit_{1, 1}$ & $\begin{pmatrix}
			\ta+2pq & 0 & 0 & \ta \\ 
			0 & \tb & \ta & 0 \\ 
			0 & \ta & \tb & 0 \\ 
			\ta & 0 & 0 & \ta-2pq
			\end{pmatrix}$ & Table \ref{tab:h11} & \S \ref{subsec:hit11} & $2p^2 = -1 $ or $ 2q^2 = 1$
			\\ \hline
			$\hit_{1, 2}$ & $\begin{pmatrix}
			p & 0 & 0 & k \\ 
			0 & p & p-q & 0 \\ 
			0 & 0 & q & 0 \\ 
			0 & 0 & 0 & -q
			\end{pmatrix}$ & Tables \ref{tab:h12}, \ref{tab:h12a},\newline \ref{tab:h12c}, \ref{tab:h12ac} 
			& \S \ref{subsec:hit12} & $p = -1 $ or $ q = 1 $ 
			\\ \hline
			$\hit_{1, 3}$ & $\begin{pmatrix}
			k^2 & kp & -kp & pq \\ 
			0 & k^2 & 0 & kq \\ 
			0 & 0 & k^2 & -kq \\ 
			0 & 0 & 0 & k^2
			\end{pmatrix}$ & Table \ref{tab:h13} & \S \ref{subsec:hit13} & $ k^4=1 $
			\\ \hline
			$\hit_{1, 4}$ & $\begin{pmatrix}
			0 & 0 & 0 & p \\ 
			0 & 0 & k & 0 \\ 
			0 & k & 0 & 0 \\ 
			q & 0 & 0 & 0
			\end{pmatrix}$ 
			& Table \ref{tab:h14}& \S \ref{subsec:hit14} & $k=-1$ or $ pq=1 $ \\\hline 
			$\hit_{0, 1}$ & $\begin{pmatrix}
			k & 0 & 0 & k \\ 
			0 & -k & 0 & 0 \\ 
			0 & 0 & -k & 0 \\ 
			0 & 0 & 0 & k
			\end{pmatrix}$ & \eqref{eqn:Nichols10,k1},  \eqref{eqn:Nichols10,k-1} & \S \ref{subsec:hit01} & $ k^2=1 $ 
			\\ \hline
			\end{tabular}
			\end{center}
			\end{table}

Here is our main result:

\begin{theorem}\label{th:main}
Let $(V,c)$ be a rigid braided vector space of dimension 2, not of diagonal type. 
If $\toba(V)$ has quadratic relations, then $(V,c)$ is as in Table \ref{tab:general}.
Furthermore, the explicit presentation of $ \toba(V)$, a PBW-basis, the dimension and the GK-dimension are given.
\end{theorem}

In Table \ref{tab:general}, the first column follows the conventions of \cite{H}, cf. \S \ref{subsec:intro_with_Hietarinta};
the third provides the defining relations, while the fourth contains the reference to the proof
and the fifth, the conditions needed to have quadratic relations.  We use for $\hit_{1,1}$ the notation
\begin{align}\label{eqn:taandtb}
\ta &= p^2-q^2, & \tb &= p^2+q^2.
\end{align}

\subsection{} \label{sec:intro_discussion} 
We discuss some features of the outcome of Theorem \ref{th:main}. 
First, in all cases where at least one non-trivial quadratic relation exists,
the Nichols algebra can be determined. This is neatly different from the diagonal case.
Indeed, consider a braided vector space $V$ of diagonal type with generalized Dynkin diagram
$\xymatrix{\circ^{q} \ar@{-}[r]^{p} & 	\circ^{-1} }$ (see \cite{H-classif RS} for details). Then
$\toba(V)$ has the quadratic relation $x_2^2 =0$ but to the best of our knowledge 
the complete set of defining relations is in general not known.

Second, the list of (isomorphism classes of) algebras underlying the Nichols algebras in our classification
is rather short, see the Appendix, and the same algebra underlies various different Nichols algebras, as
already present in the motivating examples from \cite{GGi}. This suggests relations between the different Nichols algebras,
e.g. twisting.

\subsection{} \label{sec:intro_aplication} 
Let $(V,c)$ be a rigid braided vector space and let $R = \tau c$ be the corresponding solution of the QYBE, cf. 
\S \ref{subsec:intro_with_Hietarinta}.
For applications of the determination of $\toba(V)$, 
it is necessary to find Hopf algebras $H$ that realize $V$ in $\ydh$, see \S \ref{subsec:ydh}.
Indeed, $\toba(V)$ becomes a Hopf algebra in $\ydh$ and the bosonization $\toba(V) \# H$
is a Hopf algebra that likely inherits  properties of $\toba(V)$ and $H$. 
For instance, under some mild conditions, cf. \cite[Lemma 2.2]{AAH},
\begin{align*}
\GK \toba(V)\# H = \GK \toba(V) + \GK H.
\end{align*} 
Notice that $\GK \toba(V) \leq 2$ for all $V$ in Table \ref{tab:general}.

Now, if $H$ is a co-quasitriangular Hopf algebra, then the category $\mathcal M^H$
of right $H$-comodules embed into $\ydh$ \cite[(3.4)]{T-survey}.
Thus we may seek for co-quasitriangular Hopf algebras $H$ that realize $(V, c)$ in $\mathcal M^H$.
The celebrated FRT-construction \cite{FRT} provides a universal co-quasitriangular \emph{bialgebra} $A(R)$ 
realizing $(V, c)$; see \cite{T-survey} for details and more references.
When $\dim V = 2$, $A(R)$ was computed explicitly in \cite{ACDM1,ACDM2,FHR,T}.

Further, the rigidity of $c$ is equivalent to the existence of a Hopf algebra $H(R)$ realizing 
$(V,c)$ as a $H(R)$-comodule \cite{Sch,Ha}. 
The construction of $H(R)$ is simpler when $A(R)$ has a group-like element $g$ with favourable properties
(the so-called quantum determinant)
so that $H(R)$ is the localization of $A(R)$ at $g$ \cite{FRT,T}; but otherwise is involved. 
Particularly, it is not evident to us whether  $\GK H(R) < \infty$ if $\GK A(R) < \infty$, except when localizing at the
quantum determinant.
When $\dim V = 2$,  $H(R)$ is a localization of $A(R)$ by the quantum determinant for
$ \hit_{2, 1}$ \cite{T}, but otherwise it is not known, apparently.

\subsection*{Notation}
If $j \le k \in \N_0$, then we denote $\I_{j, k} =\{j, j+1, \dots, k\}$ and $\I_k = \I_{1,k}$.
We denote by $\G_N$, $N\geq 1$, the group of $N$-th roots of unity and by $\G_N'\subset \G_N$ the subset of primitive roots. 
We also set $ \G_\infty = \bigcup_{N\geq 1} \G_N$, $ \G'_\infty = \G_\infty - \{1\}$.
Given $ n\in\mathbb{N} $, $ 0\leq i \leq n $ and $ q\in\k $, we consider the $ q $--numbers
\begin{align*}
(n)_q &= \sum_{j=0}^{n-1} q^j, &(n)_q^! &= \prod_{j=1}^{n} (j)_q, &\binom{n}{i}_q &= \frac{(n)_q^!}{(n-i)_q^!(i)_q^!}
\end{align*}

The multiplication of an algebra $A$ is denoted $\mu$, or $\mu_A$ when emphasis is needed.

\subsection*{Acknowledgments} We thank   Iv\'an Angiono,
Juan Cuadra and Leandro Vendramin for fruitful conversations at different stages of our research.

\section{Preliminaries}\label{sec:preliminaries}
\subsection{Yetter--Drinfeld modules} \label{subsec:ydh}
For more information on the material in this Subsection, see e.g. \cite{AHS, AS2}. 
Let $H$ be a Hopf algebra. 
A Yetter-Drinfeld module $V$ over $H$ is simultaneously a left $H$-module and a left $H$-comodule via $\rho$ satisfying the  compatibility 
\begin{align*}
\rho(h\cdot v) &= h_{(1)}v_{(-1)}\mathcal{S}(h_{(3)})\ot h_{(2)}\cdot v_{(0)}, &v&\in V, \, h\in H.
\end{align*}
The category of Yetter-Drinfeld modules over $H$, denoted by $\ydh$, is a braided tensor one:
the tensor product of $V, W\in \ydh$ has the natural action and coaction, while the braiding is given by
\begin{align*}
c_{V,W}(v\cdot w) &= v_{(-1)} \dot w \ot v_{(0)}, &v&\in V, \, w\in W.
\end{align*}
Thus, every $V \in \ydh$ becomes a 
rigid braided vector space $(V, c_{V, V})$; 
but a rigid braided vector space could be realized in $\ydh$ for many possible $H$. 
As every braided tensor category, $\ydh$ has some special features:
\begin{itemize}[leftmargin=*]\renewcommand{\labelitemi}{\tiny $\diamondsuit$}
	\item Algebras in $\ydh$ are objects in $\ydh$ that are also associative algebras, with the unit and multiplication maps being morphisms in $\ydh$. Same for coalgebras.
	
	\item The tensor product of two algebras $A$ and $B$ in $\ydh$ is the object $A \ot B \in \ydh$ with the multiplication $\mu_{A \ot B} = (\mu_{A} \ot \mu_{B}) (\id_A \ot c_{B, A} \ot \id_B)$. It is denoted $A \underline{\ot} B$.
	
	\item A bialgebra in $\ydh$ is an algebra and coalgebra $R$ in $\ydh$ with the comultiplication $\Delta: R \to R \underline{\ot} R$ and the counit being morphisms of algebras. 
	
\item 	A Hopf algebra in $\ydh$ is a bialgebra whose identity is convolution-invertible (the inverse is called the antipode).
\end{itemize}

Here is a distinguished fact of $\ydh$: if $R$ is a Hopf algebra in $\ydh$, then $R \#H =$ the vector space $R \ot H$ 
with the semi-direct product and semi-direct coproduct structures is again a Hopf algebra called the (Radford-Majid) bosonization of $R$ by $H$.

\subsection{Nichols algebras} \label{subsec:nichols}
Let $V \in \ydh$. The tensor algebra $T(V)$ has a natural structure of (graded) Hopf algebra in $\ydh$,
where the comultiplication is the algebra map $\Delta: T(V) \to T(V) \underline{\ot} T(V)$ given by 
$\Delta(v) = v \otimes 1 + 1 \ot v$, $v\in V$.
Let $\Cg$ be the set of Hopf ideals of $T(V)$ that are generated by homogeneous elements of degree $\geq 2$ and subobjects in $\ydh$.

\begin{definition}\label{def:nichols}
	The Nichols algebra $\toba(V)$ is the quotient of the $T(V)$ by the largest ideal $\J(V)$  in $\Cg$.
\end{definition}

Thus, if $I\in \Cg$, then $\bB = T(V)/I$ is a graded Hopf algebra in $\ydh$ 
and there is a natural epimorphism $\bB \to \toba(V)$; such $\bB$ are called pre-Nichols algebras over $V$.

We emphasize $\toba(V, c)$, $\J(V, c)$, if needed.
As a coalgebra, $\toba(V)$ depends only on the braided vector space $(V,c)$ and not on the realization in $\ydh$. 
The ideal $\J(V) = \oplus_{n \geq 2} \J^n(V)$ has alternative descriptions. 
For instance, $\J^n(V)$ is the kernel of the quantum symmetrizer associated to $c$; in particular
\begin{align}\label{eq:J2}
\J^2(V) &= \ker (\id + c).
\end{align}

Left skew derivations are also quite useful. 
Let $f \in V^*$ and consider $\partial^L_{f} = \partial_{f} \in \End  T(V)$ given by
\begin{align}\label{eq:derivada}
\begin{aligned}
	\partial_f(1) &=0,& &\partial_f(v) = f(v), \forall v\in V,
\\
\partial_f(xy) &= \partial_f(x)y + \sum_i x_i \partial_{f_i}(y), & &\text{where } c^{-1} (f \ot x) = \sum_i x_i \ot f_i.
\end{aligned}\end{align}

Let $\bB$ be a pre-Nichols algebra over $V$. Then for  every $f \in V^*$,
there is $\partial_{f} \in \End  \bB$ satisfying \eqref{eq:derivada}, compatible with the projection
$T(V) \to \bB$. This was proved in \cite{AAH} albeit previous partial formulations appeared much before.

We now outline the method  to determine Nichols algebras, already used in many papers.

First, one looks for relations in $\J(V)$: the quadratic ones are given by \eqref{eq:J2}; if $c(v \ot v) = q v \ot v$ for some 
$q\in \G'_n$, $N \ge2$, then $v^N \in \J(V)$; for more relations the following fact is convenient:
Let $x \in T(V)$. Then
\begin{align}\label{eq:deriv-criteria} 
\partial_{f} (x) =0  \text{ for all } f\in V^* \implies x \in \J(V).
\end{align}

So, one disposes of a collection of homogeneous relations $(r_\ell)_{\ell \in \Lambda}$. Let $I = \langle r_\ell: \ell \in \Lambda\rangle$. 
We may assume that $I$ is a Hopf ideal, since usually the $r_{\ell}$'s are either primitive, or primitive modulo other relations.

Second, one considers the pre-Nichols algebra $\bB = T(V)/ I$ and the natural projection $\pi: \bB \twoheadrightarrow \toba(V)$.
Assume that there is $0 \neq r\in \ker \pi$ that we choose homogeneous of the lowest degree. Then $\partial_{f} (r) \in \ker \pi$,
hence $\partial_{f} (r) =0$ by minimality of the degree. If this leads to a contradiction, then $\pi$ is actually an isomorphism.

\subsection{The braid equation in rank 2} \label{subsec:intro_with_Hietarinta}

We say that $R \in GL(V \otimes V)$  satisfies the quantum Yang-Baxter equation (QYBE)
\begin{align}\label{eqn:qybe}
R_{12}R_{13}R_{23} &=  R_{23}R_{13}R_{12}.
\end{align}
Let $\tau: V\ot V \to V \ot V$ be the usual flip, $\tau(x \ot y) = y \ot x$.

The QYBE \eqref{eqn:qybe} is equivalent to the braid equation \eqref{eqn:braid-bis} in two ways:
\begin{align}
R \overset{\star}{\longleftrightarrow} c &= \tau R, & R \overset{\ast}{\longleftrightarrow} \widetilde{c} &= R\tau. 
\end{align}

Therefore, we have involutions of the set of solutions of the QYBE \eqref{eqn:qybe}, respectively the braid equation \eqref{eqn:braid-bis}, given by
\begin{align}
R \mapsto  R^{\#} &:= \tau R \tau, & c \mapsto  c^{\#} &:= \tau c \tau. 
\end{align}
Obviously, $c = \tau R$ iff $c = R^{\#}\tau$, etc.

\subsubsection{Classification} \label{subsubsec:classif-Hietarinta}
Assume that $\dim V =2$.
The classification of the solutions of \eqref{eqn:qybe}  was performed by J. Hietarinta  in \cite{H},
up to the equivalence generated by the following relations, listed exactly as in \emph{loc. cit.}:

\begin{enumerate}[leftmargin=*,label=\rm{(\alph*)}] 
	\item\label{item:hiet-a} $R \leftrightarrow R^t$ (the transpose). 
	
\smallbreak
	\item\label{item:hiet-b}  The change of basis $R \mapsto (F\ot F) R (F\ot F)^{-1}$, $F\in GL(V)$ given by 
	$x_1 \mapsto x_2$, $x_2 \mapsto x_1$; a particular case of \ref{item:hiet-e} below.
	
\smallbreak	\item\label{item:hiet-c}  A change of basis given by
$x_1 \ot x_1 \mapsto x_1 \ot x_1$,  $x_2 \ot x_1 \mapsto x_1 \ot x_2$, 
$x_1 \ot x_2 \mapsto x_2 \ot x_1$,  $x_2 \ot x_2 \mapsto x_2 \ot x_2$.
This is just the assignment $R \mapsto R^{\#}$ defined above.	 
	
\smallbreak
\item\label{item:hiet-d} Homotheties $R \mapsto \kappa R$, $\kappa \in \ku^\times$.
	
\smallbreak
\item\label{item:hiet-e} The change of basis $R \mapsto (\varphi\ot \varphi) R (\varphi\ot \varphi)^{-1}$, $\varphi\in GL(V)$.
\end{enumerate}

 Thus, we proceed case-by-case following the classification in \cite{H} up to $\star$.
It turns out that there are 23 families of solutions, numbered $(\hit_{0, j})_{j\in \I_6}$, $(\hit_{1, j})_{j\in \I_{12}}$, $(\hit_{2, j})_{j\in \I_4}$ and $\hit_{3, 1}$.
We are interested in those invertible and rigid; we also exclude the diagonal braidings. 
The remaining solutions are listed in Table \ref{tab:general}. 
Notice that we homogeneize $\hit_{0,1}$ from \cite{H} for the discussion in \S \ref{subsubsec:nichols-Hietarinta} below.

\begin{rem} The $R$-matrices  
	$(\hit_{0, j})_{j= 4,5,6}$, $(\hit_{1, j})_{j= 5, \dots, 12}$ and $\hit_{2, 4}$ are not invertible;
	$(\hit_{0, j})_{j = 2,3}$ are invertible but not rigid;
	$\hit_{3, 1}$ is of diagonal type.
\end{rem}

We shall denote by $\Vhit_{i,j}$ the braided vector space corresponding to $R$-matrix $\hit_{i,j}$.
We have fixed a basis $x_1, x_2$ of $\Vhit_{i,j}$ and the $R$-matrix is described in the ordered basis 
$x_1 \ot x_1$, $x_2 \ot x_1$, $x_1 \ot x_2$, $x_2 \ot x_2$
of $\Vhit_{i,j} \otimes \Vhit_{i,j}$. 
In Section \ref{sec:quadratic-relations} we determine the quadratic relations of  these 8 families
and describe the Nichols algebras of them, when there are quadratic relations.

\subsubsection{Behaviour of Nichols algebras up to equivalence} \label{subsubsec:nichols-Hietarinta}

It remains to discuss the relations between the Nichols algebras under the equivalence in \cite{H}.
Let $R$, $R'$ be solutions of \eqref{eqn:qybe} on $V$ and $V'$ respectively 
(where $\dim V = \dim V' <\infty$), and $c = \tau R$, $c' = \tau R'$.
Assume that $c$, $c'$ are rigid. Notice that there are two ways of identifying $(V \ot V)^*$ with $V^* \ot V^*$, namely
\begin{align*}
\langle f \ot g, v \ot w\rangle & \overset{(I)}{=} \langle f, w\rangle  \langle g, v \rangle &
&\text{and} & \langle f \ot g, v \ot w\rangle & \overset{(II)}{=} \langle f, v\rangle  \langle g, w \rangle,
\end{align*}
$v, w\in V$, $f,g \in V^*$. The map $(V \ot V)^* \to (V \ot V)^*$ induced by these two identifications turns out to be the change of basis in \ref{item:hiet-c}. In short, we have:

\begin{itemize}[leftmargin=*]
	\item[\ref{item:hiet-a}] Let $R' = R^t$. Then $c^t = R^t \tau$, $(c^t)^\# = \tau R^t = c'$.
	
\smallbreak
\item[\ref{item:hiet-c}] Let $R' = R^\#$. Then $c' = R \tau$ and $(c')^\# = c$.

\smallbreak
\item[\ref{item:hiet-d}] if  $R' = \kappa R$ for $\kappa\in \k$, then $c' = \kappa c$.
Observe that if $R$ belongs to the family $\hit_{i,j}$, then so does $\kappa R$; 
this is why we homogeneize $\hit_{1,1}$.
(Actually $\J^2(V, c') = \ker (\id + c') = \ker (\id + \kappa c) = \ker (\kappa^{-1} +  c)$.
Therefore, $\J^2(V, c') \neq 0 \iff -\kappa^{-1}$ is an eigenvalue of $c$).

\smallbreak
\item[\ref{item:hiet-e}] Let $\varphi\in GL(V)$.
If  $R' = (\varphi\ot \varphi) R (\varphi\ot \varphi)^{-1}$, then $c' = (\varphi\ot \varphi) c (\varphi\ot \varphi)^{-1}$. Hence the map $T(\varphi): T(V) \to T(V)$ is an algebra isomorphism that induces a Hopf algebra isomorphism
$\toba(\varphi): \toba(V,c) \to \toba(V,c')$.
\end{itemize}

In conclusion, the new braided vector spaces to be considered  from a solution $ \hit_{i, j}$
 arise by the three
 transformations \ref{item:hiet-a}, \ref{item:hiet-c} and \ref{item:hiet-a} composed with \ref{item:hiet-c}.
 Indeed the transformations \ref{item:hiet-a} and \ref{item:hiet-c} are commuting involutions.
 Sometimes we will use \ref{item:hiet-b}, when it coincides with any of those.

\section{Nichols algebras of rank two with quadratic relations}\label{sec:quadratic-relations}

This Section contains the case-by-case analysis of the Nichols algebras with quadratic relations
of the braided vector spaces $(V,c)$ in the list in \cite{H}. We fix a basis $(x_i)_{i\in \I_2}$ of $V$;
let $(f_i)_{i\in \I_2}$ be its dual basis. We denote $\partial_i = \partial_{f_i}$.

\subsection{Case $ \hit_{2, 1}$}\label{subsec:hit21}
 We assume that $ k, p, q\neq 0 $ and $ k^2\neq pq $. The associated braiding is
\begin{align*}
(c(x_i \ot x_j))_{i, j \in \I_2} = \begin{pmatrix}
k^2 x_1\otimes x_1 & kq x_2\otimes x_1 + (k^2-pq) x_1\otimes x_2 \\
kp x_1\otimes x_2 & k^2 x_2\otimes x_2
\end{pmatrix}.
\end{align*}

\begin{rem}
The braided vector space $ K^1 = \Bbbk\{z_0, z_1, f_0, f_1\} $ associated to the Nichols algebra $K$ of \cite[Lemma 4.37]{AAH} has the following braiding
\begin{gather*}
\begin{aligned}
c(z_0\ot z_0)&= -z_0\ot z_0, &
c(z_0\ot z_1)&= -q_{21}(z_1+f_0)\ot z_0,\\
c(z_0\ot f_0)&= -q_{21}f_0\ot z_0, &
c(z_0\ot f_1)&= -q_{21}^2f_1\ot z_0,
\end{aligned}
\\\begin{aligned}
c(z_1\ot z_0) &= -q_{12}z_0\ot z_1 -(2z_1+f_0)\ot z_0, &
c(z_1\ot z_1)&= -z_1\ot z_1,\\
c(z_1\ot f_0)&= -f_0\ot z_1+2q_{21}f_1\ot z_0, &
c(z_1\ot f_1)&= q_{21}f_1\ot z_1,\\
c(f_0\ot z_0) &= -q_{12}z_0 \ot f_0-2f_0\ot z_0, &
c(f_0\ot f_0)&= -f_0\ot f_0, \\
c(f_0\ot z_1)&= -z_1\ot f_0-2q_{21}f_1\ot z_0, &
c(f_0\ot f_1)&= q_{21}f_1\ot f_0,
\end{aligned} \\ \begin{aligned}
c(f_1\ot z_0)&= -q_{12}^2z_0\ot f_1-2q_{12}f_0\ot z_1+2q_{12}z_1\ot f_0,\\
c(f_1\ot z_1)&=q_{12}(z_1-f_0)\ot f_1+f_1\ot (f_0-2z_1),\end{aligned} \\
\begin{aligned}c(f_1\ot f_0)&= q_{12}f_0\ot f_1-2f_1\ot f_0, &
c(f_1\ot f_1)&=-f_1\ot f_1.
\end{aligned}
\end{gather*}
Observe that $V=\Bbbk\{f_0, f_1\}$ is a braided vector subspace and $K^1/V$ is a quotient braided space
that fit in this case. 
\end{rem}

\begin{rem}\label{rem:GGi-uno}
The braidings of the braided vector spaces $ V_{i,0} $ and $ V_{i,2} $ considered in \cite[Prop. 3.4]{GGi} belong to this case.
\end{rem}

Here is our main result for this family.
Let $N = \begin{cases} \ord k^2, &\text{if } 1 \neq \ord k^2; \\ \infty,  &\text{otherwise.}
\end{cases}$

\begin{pro} \label{prop:case2,1}
If $ k^2 \neq -1 $ and $ pq \neq 1$, then there are no quadratic relations. Otherwise, the Nichols algebras are
\begin{align}\label{eqn:Nichols2,1}
&&\toba(V) &=T(V)/ \lg  x_1x_2 - kqx_2x_1, r_1, r_2\rg &&
\\\label{eqn:prop2.1} 
&\text{ \hspace{-1.3cm} where }&  r_i&= x_{i}^{N},\quad  i =1, 2  \textrm{ only if }  N < \infty.
\end{align}
Also $B = \{x_2^{a_2}x_1^{a_1}: \, 0\leq a_i < N\}$  is a PBW-basis of $ \toba(V) $ and 
\begin{align*}
\GK \toba(V) = \begin{cases} 0 &\text{if } N < \infty\ (\text{and } \dim \toba(V) = N^2), \\ 2  &\text{otherwise.}\end{cases}
\end{align*} 
\end{pro}

\pf Set $u = \lambda_1x_1^2+\lambda_2x_1x_2+\lambda_3x_2x_1+\lambda_4x_2^2$. Then
\begin{align}\label{eqn:2.1}
\Delta(u) &= u\ot 1 + 1\ot u + \lambda_1(1+k^2)x_1\ot x_1 + \lambda_4(1+k^2)x_2\ot x_2 \\
\notag &+ (\lambda_3+\lambda_2kq)x_2\ot x_1 + (\lambda_3kp+\lambda_2(1+k^2-pq))x_1\ot x_2.
\end{align}

From \eqref{eqn:2.1},  all assertions made about quadratic relations hold. The relations \eqref{eqn:prop2.1} are also clear. 
Let $\bB$ be the pre-Nichols algebra in the right-hand side of \eqref{eqn:Nichols2,1} and let $\pi: \bB\to \toba(V)$
be the natural projection. A standard argument shows that $B$ is a system of linear generators of $\bB$. 
We claim that the image of $B$ (that we call $B$ again) in $\toba(V)$ is linearly independent; 
the claim implies that $\pi$ is an isomorphism and that $B$ is a basis as claimed.
Assume that $B$ is not linearly independent and pick $r$ a homogeneous relation of minimal degree $n > 2$:
\begin{align*}
r = \sum_{M_2\le a \le M_1} c_a x_2^ax_1^{n-a},
\end{align*}
where $M_1 = \min \{n, N-1\}$, $M_2 = \max \{0, n + 1 - N\}$. Then
\begin{align*}
\Delta(r)&= \sum_{M_2\le a \le M_1} c_a \Big(\sum_{i=0}^a \binom{a}{i}_{k^2} x_2^i\ot x_2^{a-i}\Big)
\Big(\sum_{j=0}^{n-a} \binom{n-a}{j}_{k^2} x_1^j\ot x_1^{n-a-j}\Big) \\
&= \sum_{M_2\le a \le M_1}\sum_{i=0}^a\sum_{j=0}^{n-a} c_a \binom{a}{i}_{k^2}\binom{n-a}{j}_{k^2} (kp)^{(a-i)j} x_2^ix_1^j\ot x_2^{a-i}x_1^{n-a-j}.
\end{align*}
Thus
\begin{align*}
\partial_1(r) &= \sum_{M_2\le a \le \min \{M_1, n-1\}} c_a (n-a)_{k^2} (kp)^{a}\, x_2^{a}x_1^{n-a-1},
\\ \partial_2(r) &= \sum_{\max\{M_2, 1\} \le a \le M_1} c_a (a)_{k^2}  x_2^{a-1}x_1^{n-a}.
\end{align*}
By minimality of $n$, $\partial_1(r) = \partial_2(r) =0$, hence
\begin{align*}
c_a (n-a)_{k^2} &= 0,& M_2 &\le a \le \min \{M_1, n-1\}, \\ c_a (a)_{k^2} &= 0,& \max\{M_2, 1\} &\le a \le M_1.
\end{align*}
If any $ c_a\neq 0 $, then either $(a)_{k^2}$ or $(n-a)_{k^2} = 0$, but this contradicts 
the definition of $N$. Therefore the coefficients $ c_a $ are trivial and $\pi$ is bijective.
\epf

We close this Subsection by a discussion of the Nichols algebras arising from the equivalence in \cite{H}. 
First, \ref{item:hiet-b} and \ref{item:hiet-c} give rise to the same braiding, that is
\begin{align}\label{eq:braiding21bc}
(c'(x_i \ot x_j))_{i, j \in \I_2} = \begin{pmatrix}
k^2 x_1\otimes x_1 & kp x_2\otimes x_1 \\
kq x_1\otimes x_2 + (k^2-pq) x_2\otimes x_1  & k^2 x_2\otimes x_2
\end{pmatrix}.
\end{align}
Since  \ref{item:hiet-b} is a change of basis, the Nichols algebras are isomorphic.
Second,  \ref{item:hiet-a} gives rise to the braiding
\begin{align}\label{eq:braiding21a}
(c''(x_i \ot x_j))_{i, j \in \I_2} = \begin{pmatrix}
k^2 x_1\otimes x_1 & kq x_2\otimes x_1 \\
kp x_1\otimes x_2 + (k^2-pq) x_2\otimes x_1  & k^2 x_2\otimes x_2
\end{pmatrix}.
\end{align}
But \eqref{eq:braiding21a} is \eqref{eq:braiding21bc} up to $p \leftrightarrow q$, so no new Nichols algebra arises.

Third, \ref{item:hiet-a} composed with \ref{item:hiet-c} gives the initial $ \hit_{2, 1}$  
up to $p \leftrightarrow q$, so no new Nichols algebra arises.

\subsection{Case $ \hit_{2, 2}$}\label{subsec:hit22} 
We assume that $ k, p, q\neq 0 $ and $ k^2\neq pq $. The associated braiding is
\begin{align*}
(c(x_i \ot x_j))_{i, j \in \I_2} = \begin{pmatrix}
k^2 x_1\otimes x_1 & kq x_2\otimes x_1 + (k^2-pq) x_1\otimes x_2 \\
kp x_1\otimes x_2 & -pq x_2\otimes x_2 
\end{pmatrix}.
\end{align*}
Let $N_1 = \begin{cases} \ord k^2, &\text{if } 1 \neq \ord k^2; \\ \infty,  &\text{otherwise.}
\end{cases}$ and $N_2 = \begin{cases} \ord (-pq), &\text{if } 1 \neq \ord (-pq); \\ \infty,  &\text{otherwise.}
\end{cases}$

\begin{pro} If $ k^2 \neq -1 $ and $ pq \neq 1$, then there are no quadratic relations.
	Otherwise, the Nichols algebras are
\begin{align}\label{eqn:Nichols2,2}
\toba(V) &= T(V) / \lg  x_1x_2 - kqx_2x_1, r_1, r_2\rg 
\end{align}
where $ r_i= x_i^{N_i}, \, i\in \I_2$, only if $ N_i<\infty $. Also $ \{ x_2^{a_2}x_1^{a_1}:  a_i \in \I_{0, N_i-1}\}$ 
is a PBW-basis and $ \GK \toba(V) = |\{i\in\I_2: N_i =\infty\}|$; if 0, then $\dim \toba(V) = N_1N_2$.
\end{pro}

\pf Similar to the proof of Proposition \ref{prop:case2,1}.
\epf

We next discuss the Nichols algebras arising from the equivalence in \cite{H}. 
 First,  \ref{item:hiet-a} gives rise to the braiding
\begin{align}\label{eq:braiding22a}
(c'(x_i \ot x_j))_{i, j \in \I_2} = \begin{pmatrix}
k^2 x_1\otimes x_1 & kq x_2\otimes x_1  \\
kp x_1\otimes x_2 + (k^2-pq) x_2\otimes x_1 & -pq x_2\otimes x_2 
\end{pmatrix}.
\end{align}

Let $N_1$ and $N_2$ be as above.

\begin{pro} If $ k^2 \neq -1 $ and $ pq \neq 1$, then there are no quadratic relations.
	Otherwise, the Nichols algebras are
	\begin{align}\label{eqn:Nichols2,2a}
	\toba(V) &=T(V)/ \lg  x_2x_1 - kpx_1x_2, r_1, r_2\rg 
	\end{align}
where $ r_i= x_i^{N_i}, \, i\in \I_2$, only if $ N_i<\infty $. Also $ \{ x_2^{a_2}x_1^{a_1}:  a_i \in \I_{0, N_i-1}\}$ 
is a PBW-basis and $ \GK \toba(V) = |\{i\in\I_2: N_i =\infty\}|$; if 0, then $\dim \toba(V) = N_1N_2$.
\end{pro}

\pf Similar to the proof of Proposition \ref{prop:case2,1}.
\epf

Second, \ref{item:hiet-c} gives rise to the braiding
\begin{align}\label{eq:braiding22c}
(c''(x_i \ot x_j))_{i, j \in \I_2} = \begin{pmatrix}
k^2 x_1\otimes x_1 & kp x_2\otimes x_1  \\
kq x_1\otimes x_2 + (k^2-pq) x_2\otimes x_1 & -pq x_2\otimes x_2 
\end{pmatrix}.
\end{align}
But \eqref{eq:braiding22c} is \eqref{eq:braiding22a} up to $p \leftrightarrow q$, so no new Nichols algebra arises here.

Third, \ref{item:hiet-a} composed with \ref{item:hiet-c} gives the initial $ \hit_{2, 1}$  
up to $p \leftrightarrow q$, so no new Nichols algebra arises.

\subsection{Case $ \hit_{2, 3}$}\label{subsec:hit23} We assume that $ k\neq 0$, and either  $p\neq 0$, or $q\neq 0$, or $s\neq 0$. 
The associated braiding is $(c(x_i \ot x_j))_{i, j \in \I_2} =$
\begin{align*}
= \begin{pmatrix}
k x_1\otimes x_1 & k x_2\otimes x_1 + q x_1\otimes x_1 \\
k x_1\otimes x_2 + p x_1\ot x_1 & k x_2\otimes x_2 + s x_1\otimes x_1 + p x_2\otimes x_1 + q x_1\otimes x_2
\end{pmatrix}.
\end{align*}

\begin{rem}
The braided vector spaces $\mathcal{V}(\eps, 2)$ considered in \cite[\S 1.2]{AAH} 
fit in this case taking $k=\eps$, $q=1$ and $p = s = 0$. 
In particular, the Jordan plane $\mathcal{V}(1, 2)$ and the super Jordan plane 
$\mathcal{V}(-1, 2)$ belong to this case. 
\end{rem}

To state the next result, we need the notation 
\begin{align*}
x_{21} = [x_2, x_1]_c = x_2x_1 - \mu (c(x_2\ot x_1)).
\end{align*}

\begin{table}[ht]
	\caption{Nichols algebras of type $ \hit_{2, 3} $}\label{tab:h23}
	\begin{center}
		\begin{tabular}{| c | p{1cm} | c | p{4,9cm} | c | c |}\hline
			$k$ & $p$ & $ s $ & $ \J(V)$ & Basis &$\GK$
			\\\hline
			\begin{small}
				$-1$ \end{small} & $ -q $ & $ q^2 $ &$ \lg x_1^2, x_2^2 - qx_1x_2, x_1x_2 + x_2x_1 \rg$ & $ (*_{2, 3})_1 $ 
			& \begin{small} $0, \dim = 4$ \end{small}
			\\\cline{3-6}
			 &  & $ \neq q^2 $ &$\lg x_1^2, x_1x_2 + x_2x_1 \rg$ & $ (*_{2, 3})_2 $ & $1$
			\\\cline{2-6}
			 & \begin{small}  $ \neq -q $\end{small} &  & \begin{small} 
			 	$\lg x_1^2, x_2x_{21} + (p-q)x_1x_{21} - x_{21}x_2 \rg$\end{small} &  $ (*_{2, 3})_3 $& $2$
			\\\hline
			$1$ &  &  &$\lg \dfrac{q-p}{2}x_1^2 - x_1x_2 +x_2x_1 \rg$ & $ (*_{2, 3})_4 $ & $2$
			\\\hline 
		\end{tabular}
	\end{center}
\end{table}

\begin{pro} \label{prop:Nichols23}
If $k \neq \pm 1$, then there are no quadratic relations. 
Otherwise, the Nichols algebras are as in Table \ref{tab:h23}, where
\begin{align*}
 (*_{2, 3})_1 &=  \{x_1^{a_1}x_2^{a_2}: 0\leq a_i \leq 1 \} ;\\ 
(*_{2, 3})_2 &=  \{x_1^{a_1}x_2^{a_2}: 0\leq a_1 \leq 1, 0\leq a_2 < \infty \};\\
(*_{2, 3})_3 &=  \{x_1^{a}x_{21}^{b}x_2^{c}: 0\leq a \leq 1, 0\leq b, c < \infty \}; \\
(*_{2, 3})_4 &= \{x_1^{a_1}x_2^{a_2}: 0\leq a_i < \infty \}.
\end{align*}
\end{pro}

\pf Set $u = \lambda_1x_1^2+\lambda_2x_1x_2+\lambda_3x_2x_1+\lambda_4x_2^2$. Then
\begin{align*}
\Delta(u) &= u\ot 1 + 1\ot u + (\lambda_1(1+k)+\lambda_2q+\lambda_3p+\lambda_4s)x_1\ot x_1 \\
&+ \lambda_4(1+k)x_2\ot x_2 + (\lambda_2k+\lambda_3+\lambda_4p)x_2\ot x_1 \\
&+ (\lambda_2+\lambda_3k+\lambda_4q)x_1\ot x_2.
\end{align*}
Then all the assertions on quadratic relations hold. The claim in the 
first row of Table \ref{tab:h23} follows then easily. 
To simplify the discussion of the rest, we consider three cases:
\begin{enumerate}[leftmargin=*,label=\rm{(\roman*)}] 
\item \label{case1:23} $k=-1$, $ p=-q $ and $ s\neq q^2 $;
\item \label{case2:23} $k=-1$, $ p\neq -q $;
\item \label{case3:23} $k=1$.
\end{enumerate}
Case \ref{case1:23}: We first prove by induction that for $n\geq 2$
\begin{align*}
\partial_1(x_2^n) &= \begin{cases}
\dfrac{n}{2}qx_2^{n-1}+(\dfrac{n(n-2)}{4}q^2 +\dfrac{n}{2}s)x_1x_2^{n-2}, & \textrm{if } n \textrm{ is even};\\
-\dfrac{n-1}{2}qx_2^{n-1}+ \Big(\dfrac{(n-1)(n-3)}{4}q^2  & \hspace{-20pt}+\dfrac{n-1}{2}s\Big)x_1x_2^{n-2}, 
\\& \textrm{if } n \textrm{ is odd}.
\end{cases} \\
\partial_2(x_2^n)&= \begin{cases}
-\dfrac{n}{2}q x_1x_2^{n-2}, & \textrm{if } n \textrm{ is even};\\
x_2^{n-1}-\dfrac{n-1}{2}qx_1x_2^{n-2}, & \textrm{if } n \textrm{ is odd}.
\end{cases}
\end{align*}
From this, we prove again by induction that,  for $n\geq 2$
\begin{align*}
\partial_1(x_1x_2^{n-1})&= \begin{cases}
x_2^{n-1}+\dfrac{n}{2}qx_1x_2^{n-2}, & \textrm{if } n \textrm{ is even};\\
x_2^{n-1}-\dfrac{n-1}{2}qx_1x_2^{n-2}, & \textrm{if } n \textrm{ is odd}.
\end{cases}\\
\partial_2(x_1x_2^{n-1})&= \begin{cases}
-x_1x_2^{n-2} , & \textrm{if } n \textrm{ is even};\\
0, & \textrm{if } n \textrm{ is odd}.
\end{cases}
\end{align*}
Let  $\bB = T(V)/ \langle x_1^2, x_1x_2 + x_2x_1\rangle$  and let $\pi: \bB\to \toba(V)$
be the natural projection. A standard argument shows that $(*_{2, 3})_2$ generates linearly $\bB$. 
Assume that the image $B = \pi(B)$ is not linearly independent
and pick  a linear homogeneous  
relation of minimal degree $n > 2$
\begin{align*}
r = c_1x_2^n + c_2x_1x_2^{n-1}.
\end{align*}

Assume that $n$ is odd. Then
\begin{align*}
0 &= \partial_2(r) = c_1 \left(x_2^{n-1}-\dfrac{n-1}{2}qx_1x_2^{n-2}\right)
\implies c_1=0 \implies \\
0 &= \partial_1(r) = c_2\left(x_2^{n-1}-\dfrac{n-1}{2}qx_1x_2^{n-2}\right) \implies r = 0.
\end{align*}
Assume that $n$ is even. Then
\begin{align*}
0 = \partial_1(r) &= c_1\left(\dfrac{n}{2}qx_2^{n-1}+\left(\dfrac{n(n-2)}{4}q^2 +\dfrac{n}{2}s\right) x_1x_2^{n-2}\right)\\
&+c_2\left(x_2^{n-1}+\dfrac{n}{2}qx_1x_2^{n-2}\right)\\
&= \left(\dfrac{n}{2}qc_1+c_2\right)x_2^{n-1} + \left(\left(\dfrac{n(n-2)}{4}q^2 +\dfrac{n}{2}s\right)c_1+\dfrac{n}{2}qc_2\right)x_1x_2^{n-2};\\
0 = \partial_2(r) &= -c_1\dfrac{n}{2}q x_1x_2^{n-2}-c_2 x_1x_2^{n-2} =(-\dfrac{n}{2}qc_1 -c_2)x_1x_2^{n-2}.
\end{align*}
Hence
\begin{align*}
\dfrac{nq}{2}c_1+c_2 &= 0, &
\left(\dfrac{n(n-2)q^2}{4} +\dfrac{ns}{2}\right) c_1 + \dfrac{nq}{2}c_2 &= 0.
\end{align*}

The system above has non trivial solution iff $q^2 = s$. The claim in row 2 of Table \ref{tab:h23} is established. 

\smallbreak

Case \ref{case2:23}: By analogy with the super Jordan plane $\mathcal{V}(1, 2)$, 
we look for cubic relations and obtain the following one by \eqref{eq:deriv-criteria}:
\begin{align*}
0=x_2^2x_1 + (p-q)x_1x_2x_1 - x_1x_2^2 = x_2x_{21} + (p-q)x_1x_{21} - x_{21}x_2.
\end{align*}
Let  $\bB = T(V)/ \langle  x_1^2, x_2x_{21} + (p-q)x_1x_{21} - x_{21}x_2\rangle$. 
Observe that $x_1x_{21}= x_{21}x_1$ in $\bB$. 
Arguing as in \cite{AAH}, that is using the commutation relations, 
we see  that $(*_{2, 3})_3$ is a system of linear generators of $\bB$. 
We need the formulae of the derivations on the elements of this basis.
First,
\begin{align*}
\partial_1(x_{21}^b)&=b(p+q)x_1x_{21}^{b-1}, &\partial_2(x_{21}^b)&=0,& b&\geq 1.
\end{align*}
For $i\in \I_2$, $c\geq 0$, set
\begin{align*}
\partial_i(x_{2}^c) &=\partial_{i,c}=\partial_{i,c,0}+x_1\partial_{i,c,1},&
\text{where } \partial_{i,c,j} \in \ku \{x_{21}^bx_2^d: b,d \ge 0\}.
\end{align*}

Straightforward calculations show that, for $ b\geq 1$, $c\geq 0$,
\begin{align}
\partial_1(x_{21}^bx_2^c)&=x_{21}^b(\partial_{1,c}-2bq\partial_{2,c})+b(p+q)x_1x_{21}^{b-1}x_2^c, \\
\partial_2(x_{21}^bx_2^c)&=x_{21}^{b}\partial_{2,c},\\
\partial_1(x_1x_{21}^bx_2^c)&=x_{21}^bx_2^c-x_1x_{21}^b\partial_{1,c,0}+(2b+1)qx_1x_{21}^b\partial_{2,c,0}, \\
\partial_2(x_1x_{21}^bx_2^c)&=-x_1x_{21}^{b}\partial_{2,c,0},
\\
\label{eqn:superjordan}
\partial_{2,c,0} &= \begin{cases}
0, & \textrm{if } n \textrm{ is even},\\
x_2^{c-1}, & \textrm{if } n \textrm{ is odd}.
\end{cases}
\end{align}

Assume that the image of $B$ under the projection $\pi: \bB\to \toba(V)$ is not linearly independent.
Pick $r$ a non-trivial linear combination homogeneous of minimal degree   $N\geq 4$. 

Suppose first that $ N $ is odd. Then there are scalars $\lambda_{b}$, $\mu_{t}$ such that
\begin{align*}
r =  \sum_{0 \le b \le \frac{N-1}{2}} \lambda_{b}\,x_{21}^bx_{2}^{N-2b} +\sum_{0 \le t \le \frac{N-1}{2}} \mu_{t} \, x_1x_{21}^{t}x_{2}^{N-1-2t}.
\end{align*}
Applying $ \partial_2 $ to $r$, we obtain
\begin{align*}
0 &=  \sum_{0 \le b \le \frac{N-1}{2}} \lambda_{b}x_{21}^b\partial_{2,N-2b} -\sum_{0 \le t \le \frac{N-1}{2}} \mu_{t}x_1x_{21}^{t}\partial_{2,N-1-2t,0} \\
&=  \sum_{0 \le b \le \frac{N-1}{2}} \lambda_{b}x_{21}^b(\partial_{2,N-2b,0} + x_1\partial_{2,N-2b,1}) -\sum_{0 \le t \le \frac{N-1}{2}} \mu_{t}x_1x_{21}^{t}\partial_{2,N-1-2t,0}\\
&\stackrel{\eqref{eqn:superjordan}}{=}  \sum_{0 \le b \le \frac{N-1}{2}} \lambda_{b}x_{21}^bx_2^{N-1-2b} + \lambda_{b}x_1x_{21}^b\partial_{2,N-2b,1}.
\end{align*}
From this, we see that $\lambda_{b}= 0$, $b= 0, 1, \cdots, \frac{N-1}{2}$. Therefore
\begin{align*}
0 &= \partial_1(r) \stackrel{\eqref{eqn:superjordan}}{=}  \sum_{0 \le t \le \frac{N-1}{2}} \mu_{t}(x_{21}^{t}x_2^{N-1-2t}-x_1x_{21}^{t}\partial_{1,N-1-2t,0});
\end{align*}
hence $r=0$. 

Assume next that $N$ is even. Then there are scalars $\lambda_{b}$, $\mu_{t}$ such that
\begin{align*}
r =  \sum_{0 \le b  \le \frac{N}{2}} \lambda_{b}x_{21}^bx_{2}^{N-2b} +\sum_{0 \le t  \le \frac{N -2}{2}} \mu_{t}x_1x_{21}^{t}x_{2}^{N-1-2t}.
\end{align*}
Thus
\begin{align}\label{eqn:1superjordan}
\begin{split}
0 &= \partial_2(r) \stackrel{\eqref{eqn:superjordan}}{=} \sum_{0 \le b  \le \frac{N}{2}} \lambda_{b}x_1x_{21}^b\partial_{2,N-2b,1} -\sum_{0 \le t  \le \frac{N -2}{2}} \mu_{t}x_1x_{21}^{t}x_2^{N-2-2t}\\
&= \sum_{0 \le b  \le \frac{N -2}{2}} \lambda_{b}x_1x_{21}^b\partial_{2,N-2b,1} -\sum_{0 \le t  \le \frac{N -2}{2}} \mu_{t}x_1x_{21}^{t}x_2^{N-2-2t}.
\end{split}
\end{align}
Applying $\partial_1$ to $r$, we obtain for some $z \in \toba^{N-2}(V)$
\begin{align*}
0&= \sum_{0 \le b  \le \frac{N}{2}} \lambda_{b}x_{21}^b(\partial_{1,N-2b,0}-2bq\partial_{2,N-2b,0}) +\sum_{0 \le t  \le \frac{N -2}{2}} \mu_{t}x_{21}^{t}x_{2}^{N-1-2t}+x_1z\\
&\stackrel{\eqref{eqn:superjordan}}{=} \sum_{0 \le b  \le \frac{N -2}{2}} \lambda_{b}x_{21}^b\partial_{1,N-2b,0} +\sum_{0 \le t  \le \frac{N -2}{2}} \mu_{t}x_{21}^{t}x_{2}^{N-1-2t}+x_1 z.
\end{align*}
In particular,
\begin{align}\label{eqn:2superjordan}
0&= \sum_{0 \le b  \le \frac{N -2}{2}} \lambda_{b}x_{21}^b\partial_{1,N-2b,0} +\sum_{0 \le t  \le \frac{N -2}{2}} \mu_{t}x_{21}^{t}x_{2}^{N-1-2t}.
\end{align}
Also observe that if $n\geq 2$  is even, then we obtain that for some $w_i
\in \ku \{x_{21}^bx_2^d: b,d \ge 0, \, b+d =n+i-5\}$, $ i\in \I_2 $, 
\begin{align*}
\partial_{2,n,1}&= \frac{n}{2}px_{2}^{n-2} +x_{21} w_1, &\partial_{1,n,0}&= \frac{n}{2}qx_{2}^{n-1} +x_{21} w_2.
\end{align*}
Looking at the terms $ x_1x_2^{N-2} $ in \eqref{eqn:1superjordan} and $ x_2^{N-1} $
in \eqref{eqn:2superjordan}, we get 
\begin{align*}
\begin{cases}
\frac{N}{2}p\lambda_{0} - \mu_{0} =0\\
\frac{N}{2}q\lambda_{0} + \mu_{0} =0
\end{cases}
\end{align*}
whose determinant is $ \frac{N}{2}(p+q)\neq 0 $. Thus $ \lambda_{0} = \mu_{0} =0 $. Similarly, we prove that 
$\lambda_{i} = \mu_{i} =0$, $i=0, 1, \cdots, \frac{N-2}{2} $. It remains $\lambda_{\frac{N}{2}}$, but 
\begin{align*}
\partial_1(x_{21}^{\frac{N}{2}})= \frac{N(p+q)}{2}x_1x_{21}^{\frac{N-2}{2}}
\end{align*} 
hence $r=0$. 

\smallbreak
Case \ref{case3:23}: 
We start by the following claim, whose proof is straightforward:
	\begin{align}\label{eqn:1of2,3}
c(x_1^n\ot x_2)&= x_2\ot x_1^n +nq x_1\ot x_1^n,\quad n\geq 1.
	\end{align}
Let  $\bB = T(V)/ \langle \frac{q-p}{2}x_1^2 - x_1x_2 +x_2x_1\rangle$. 
By a standard argument, $(*_{2, 3})_4$ generates linearly $\bB$. 
We need the formulae of the derivations on  $(*_{2, 3})_4$.
First, we set 
\begin{align*}
\partial_2(x_2^n) &= \sum_{0\le j\le n-1} d_{j}^{(n-1)} x_1^jx_2^{n-1-j}, \quad n\geq 1,
\end{align*}
and claim that 
\begin{align}\label{eq:3of2.3}
d_{0}^{(n-1)}&=n, \quad n\geq 1.
\end{align}
Indeed, the case $ n=1 $ is clear. Assume that \eqref{eq:3of2.3} holds for $n$.
Projecting $\Delta(x_2^{n+1})$ to $V \ot \toba^{n}(V)$, we get
\begin{align*}
\sum_{i \in \I_2} x_i\ot \partial_i(x_2^{n+1})&= (x_2\ot 1)(1\ot x_2^n) + (1\ot x_2)(x_1\ot \partial_1(x_2^n) + x_2\ot \partial_2(x_2^n))\\
&= x_1\ot (x_2\partial_1(x_2^n)+px_1\partial_1(x_2^n)+sx_1\partial_2(x_2^n)+qx_2\partial_2(x_2^n)) \\
&+x_2\ot (x_2^n+x_2\partial_2(x_2^n)+px_1\partial_2(x_2^n)).
\end{align*}
Hence, by the inductive hypothesis,
\begin{align*}
\partial_2(x_2^{n+1})&= x_2^n+x_2\partial_2(x_2^n)+px_1\partial_2(x_2^n) \\
&= x_2^n+x_2(nx_2^{n-1}+\sum_{1\le j\le n-1} d_{j}^{(n-1)} x_1^jx_2^{n-1-j})+px_1\partial_2(x_2^n) \\
&= (n+1)x_2^n +x_1(x_2+\dfrac{p-q}{2}x_1)\sum_{1\le j\le n-1} d_{j}^{(n-1)} x_1^{j-1}x_2^{n-1-j} \\
&+px_1\partial_2(x_2^n),
\end{align*}
and the claim is proved.

Observe that, by \eqref{eqn:1of2,3}, for $ n\geq 3 $
\begin{align}\label{eq:1of2.3}
\partial_2(x_1^ix_2^{n-i})&=\begin{cases}
0, & \textrm{if } i=n;\\
x_1^{i}\partial_2(x_2^{n-i}), & \textrm{if } 0\leq i<n.
\end{cases}
\end{align}
Assume that the image of $B$ under the projection $\pi: \bB\to \toba(V)$ is not linearly independent. 
Pick $r = \sum_{i=0}^n c_i x_1^ix_2^{n-i}$ a non-trivial homogeneous  relation of minimal degree $ n>2$. By \eqref{eq:1of2.3}
\begin{align}\label{eq:2of2.3}
0 = \partial_2(r) = \sum_{0\le i\le n-1} c_i x_1^{i}\partial_2(x_2^{n-i}).
\end{align}
Observe that the term $ x_2^{n-1} $  appears only one time in \eqref{eq:2of2.3}. Furthermore, by \eqref{eq:3of2.3}, we can rewrite \eqref{eq:2of2.3} as
\begin{align}\label{eq:4of2.3}
0= \partial_2(r) =  c_0 nx_2^{n-1} + \sum_{1\le j\le n-1}m_jx_1^jx_2^{n-1-j} 
\end{align}
for some $ m_j\in \Bbbk $. By minimality of $n$, we obtain $ c_0 = 0 $. Similarly, we can replace \eqref{eq:4of2.3} by
\begin{align}
0= \partial_2(r) =  c_1 (n-1)x_1x_2^{n-2} + \sum_{2\le j\le n-1}m_j'x_1^jx_2^{n-1-j} 
\end{align}
what gives us $ c_1=0 $. Inductively, we get $ c_i = 0, \, 0\leq i <n $. Thus, $ r = c_n x_1^n $. But $ \partial_1(x_1^{n}) = nx_1^{n-1} $ which implies $ r = 0 $.
\epf

We next discuss the Nichols algebras arising from the equivalence in \cite{H}. 
First, \ref{item:hiet-a} and \ref{item:hiet-b} give rise to the same braiding, that is
$(c'(x_i \ot x_j))_{i, j \in \I_2} $
\begin{align*} 
= \begin{pmatrix}
k x_1\otimes x_1 + s x_2\otimes x_2 + p x_1\otimes x_2 + q x_2\otimes x_1 & k x_2\otimes x_1 + p x_2\otimes x_2 \\
k x_1\otimes x_2 + q x_2\ot x_2 & k x_2\otimes x_2 
\end{pmatrix}.
\end{align*}
Since  \ref{item:hiet-b} is a change of basis, this Nichols algebra is isomorphic to the original one.
Second, \ref{item:hiet-c} gives rise to the braiding $(c''(x_i \ot x_j))_{i, j \in \I_2} $
\begin{align*}
= \begin{pmatrix}
k x_1\otimes x_1 & k x_2\otimes x_1 + p x_1\otimes x_1 \\
k x_1\otimes x_2 + q x_1\ot x_1 & k x_2\otimes x_2 + s x_1\otimes x_1 + q x_2\otimes x_1 + p x_1\otimes x_2
\end{pmatrix}.
\end{align*}
But this is the initial $ \hit_{2, 3}$ up to $p \leftrightarrow q$, so no new Nichols algebra arises. Thus, \ref{item:hiet-a} composed with \ref{item:hiet-c} gives the braiding $c'$ above up to $p \leftrightarrow q$, so no new Nichols algebra arises.

\subsection{Case $ \hit_{1, 1}$} \label{subsec:hit11}

We assume that $ p, q\neq 0$ and $p^2\neq q^2$ with $ \ta $ and $ \tb $ as in \eqref{eqn:taandtb}. The associated braiding is $(c(x_i \ot x_j))_{i, j \in \I_2} = $
\begin{align*}
=\begin{pmatrix}
((\ta+2pq) x_1\otimes x_1 + \ta x_2\otimes x_2 & \tb x_2\otimes x_1 + \ta x_1\otimes x_2 \\
\tb x_1\otimes x_2 + \ta x_2\ot x_1 & (\ta-2pq) x_2\otimes x_2 + \ta x_1\otimes x_1
\end{pmatrix}
\end{align*}

\begin{table}[ht]
	\caption{Nichols algebras of type $ \hit_{1, 1} $}\label{tab:h11}
	\begin{center}
		\begin{tabular}{| c | c | c | c |c | c |}\hline
			$2p^2$ & $2q^2$ & $-2pq $ & $ \J(V)$ & Basis & $ \GK $
			\\\hline 
			$-1$ & $1$  &  & $ \lg  r_1, r_2, r_3 \rg $ & $ (*_{1, 1})_1 $ & $ 0$, $\dim  = 4$
			\\\hline
			$\neq -1$ & $1$  & $\in \G_N', N\geq 3$, $N\neq 4$ & $ \lg  r_1, r_2, r_{4, N} \rg$  & $ (*_{1, 1})_2 $ & $ 0$, $\dim  = 2N $
			\\\cline{3-6} 
			 &  & $\not\in \G_\infty $ &$ \lg  r_1, r_2 \rg$ & $ (*_{1, 1})_3 $ & $1$
			\\\hline 
			$-1$ & $ \neq 1$ & $\in \G_N', N\geq 3$, $N\neq 4$ & $ \lg r_1, r_3, r_{4, N} \rg$ &$ (*_{1, 1})_2 $ & $ 0$, $\dim  = 2N $
			\\\cline{3-6} 
			 &  & $\not\in \G_\infty $ & $ \lg r_1, r_3 \rg$ & $ (*_{1, 1})_3 $ & $1$
			\\\hline
		\end{tabular}
	\end{center}
\end{table}

\begin{pro}   If $ 2p^2 \neq  -1 $ and $ 2q^2 \neq  1$, then there are no quadratic relations. Otherwise, the Nichols algebras are as in Table \ref{tab:h11}, where 
\begin{align}\label{eqn:r1of11}
r_1 &= x_1^2 - \dfrac{\ta +2pq+1}{\ta}x_2^2; \\\label{eqn:r2of11}
 r_2 &= x_1x_2 - x_2x_1; \\\label{eqn:r3of11}
r_3 &= x_1x_2 +x_2x_1; \\\label{eqn:r4of11}
r_{4,N} &=\begin{cases}
x_2^{N}, & \textrm{if } N \textrm{ is odd};\\
x_1x_2^{\frac{N-2}{2}}, & \textrm{if } N \textrm{ is even}.
\end{cases}
\end{align}
\begin{align*}
&\hspace{-0.2cm}\text{and}& (*_{1, 1})_1 &=  \{1, x_1, x_2, x_2^2 \} ;\\ 
& & (*_{1, 1})_2 &= \begin{cases}
\{x_1^{a_1}x_2^{a_2}: 0\leq a_1 \leq 1, \, 0\leq a_2 < N\}, & \textrm{if } N \textrm{ is odd};\\
\{x_1x_2^{a}: 0\leq a < \frac{N-2}{2} \}\cup\{x_2^{a}: 0\leq a < \frac{N+2}{2} \}, & \textrm{if } N \textrm{ is even}.
\end{cases}\\
& & (*_{1, 1})_3 &= \{x_1^{a_1}x_2^{a_2}: 0\leq a_1 \leq 1, \, 0\leq a_2 < \infty \}.
\end{align*}
\end{pro}

\pf  To start with, observe that $ N\neq 1, 2 $ in rows $2$ and $4$ of Table \ref{tab:h11} because $\ta \neq 0$ by hypothesis. 
Moreover, $N\neq 4$ since this case is already covered by row 1.

As previously, we set $u = \lambda_1x_1^2+\lambda_2x_1x_2+\lambda_3x_2x_1+\lambda_4x_2^2$. Then
\begin{align*}
\Delta(u) &= u\ot 1 + 1\ot u + (\lambda_1(\ta+2pq+1)+\lambda_4\ta)x_1\ot x_1 \\
&+ (\lambda_1\ta+\lambda_4(\ta-2pq+1))x_2\ot x_2 + (\lambda_2(1+\ta)+\lambda_3\tb)x_2\ot x_1 \\
&+ (\lambda_2\tb+\lambda_3(1+\ta))x_1\ot x_2.
\end{align*}

Then all the assertions on quadratic relations hold. 
In particular, row $1$ is established. To simplify the exposition, we consider two cases:
\begin{enumerate}[leftmargin=*,label=\rm{(\roman*)}] 
\item \label{case1:11} $ 2p^2\neq -1 $ and $ 2q^2 = 1 $;
\item \label{case2:11} $ 2p^2= -1 $ and $ 2q^2 \neq 1 $.
\end{enumerate}
Case \ref{case1:11}: We consider the braided Hopf algebra $\widetilde{\toba} = T(V) / \lg  r_1, r_2\rg$;
clearly 
\begin{align*}
(*_{1, 1})_3 = \{x_1^{a_1}x_2^{a_2}: 0\leq a_1 \leq 1, \, 0\leq a_2 < \infty \}
\end{align*}
generates linearly $\widetilde{\toba}$.
By induction, it follows that for $n\geq 2$
\begin{align}\label{eqn:1case1of11}
\begin{aligned}
\partial_1(x_2^n)&=g_1(n)x_1x_2^{n-2}, &   \partial_2(x_2^n)&=f_1(n)x_2^{n-1},  \\
\partial_1(x_1x_2^{n-1})  &= g_2(n)x_2^{n-1},&   \partial_2(x_1x_2^{n-1}) &= f_2(n)x_1x_2^{n-2};
\end{aligned}
\end{align}
here the functions $f_j$ and $g_j$, found explicitly with \texttt{WolframAlpha}, are
\begin{align*}
f_j(n)=\dfrac{q2^np^{2n}+(-1)^j(p+q)(-2pq)^n+p}{4pq(p+q)}; \\
g_j(n)=\dfrac{q2^np^{2n}+(-1)^j(q-p)(-2pq)^n-p}{4pq(p+(-1)^{j+1}q)}.
\end{align*}
Let $0 \neq r = c_1x_2^n + c_2x_1x_2^{n-1} \in \ker (\widetilde{\toba} \to \toba(V)) $ with $ n\geq 3 $ minimal. 
Then
\begin{align*}
\partial_1(r) &= c_1g_1(n) x_1x_2^{n-2} +c_2 g_2(n) x_2^{n-1}; \\
\partial_2(r) &= c_1f_1(n) x_2^{n-1} +c_2f_2(n) x_1x_2^{n-2}.
\end{align*}
In other words, $\ker (\widetilde{\toba} \to \toba(V)) \neq 0$  
iff there exists $n \in \N$ such that either $ f_1(n) = g_1(n) = 0 $ or $ f_2(n) = g_2(n) = 0 $. But
\begin{align*}
f_j(n) = g_j(n) = 0 &\Longrightarrow -2pq\in\G_{jn}.
\end{align*}
Hence, if $-2pq\notin \G_{\infty}$, then $\widetilde{\toba} \simeq \toba(V)$ and arguing as in all preceding cases, 
$(*_{1, 1})_3$ is a basis of $\toba(V)$. This establishes the third row of Table \ref{tab:h11}.

\medbreak
Next, we assume that $-2pq\in \G_{\infty}$ and set $N = \ord (-2pq)$.
As explained above, we suppose that $ N\neq 1,2,4$. Then:
\begin{flalign}\label{eq:discussion-above}
\begin{aligned}
&\text{If $ N $ is odd, then $ f_1(N)= g_1(N)=0 $.} \\
&\text{If $m \in \N$ and $ f_1(m)= g_1(m)=0 $, then $N | m$.	} \\
&\text{If $N$ is even, then $ f_2(N/2)= g_2(N/2)=0 $. } \\
&\text{If $m \in \N$ and $ f_2(m)= g_2(m)=0 $, then $N/2 | m$.} \\
\end{aligned}		
\end{flalign}

Therefore,  the relation $r_{4,N}$, i.e. \eqref{eqn:r4of11}, holds in $\toba(V)$.
Let us consider the braided Hopf algebra 
$\toba = T(V) / \lg  r_1, r_2, r_{4,N}\rg \simeq \widetilde{\toba} / \lg  r_{4,N}\rg$;
clearly $(*_{1, 1})_2$ generates linearly $\toba$. 
Finally, $(*_{1, 1})_2$ is linearly independent in $\toba(V)$,
as any linear relation would appear in a degree higher than $\deg r_{4,N}$ by \eqref{eq:discussion-above}. 

\smallbreak
Case \ref{case2:11}: By induction, we obtain equations \eqref{eqn:1case1of11} but this time with
\begin{align*}
f_j(n)&=\dfrac{(p-q)(-2pq)^n+(-1)^jp(-2q^2)^n+(-1)^{j+1}q}{4pq(q-p)}; \\
g_j(n)&=\dfrac{(p+q)(-2pq)^n+(-1)^jp(-2q^2)^n+(-1)^{j}q}{-4pq(p+(-1)^{j}q)}.
\end{align*}
The rest of the proof is similar to the previous case.
\epf

Finally, there are no new Nichols algebras arising from the equivalence in \cite{H},
because all changes \ref{item:hiet-a}, \ref{item:hiet-c} and \ref{item:hiet-a} composed with \ref{item:hiet-c} give rise to the same initial braiding $ \hit_{1, 1}$.

\subsection{Case $ \hit_{1, 2}$}\label{subsec:hit12} 
We assume that $ p, q\neq 0$, and either  $p\neq q$, or $k\neq 0$. The associated braiding is
\begin{align*}
(c(x_i \ot x_j))_{i, j \in \I_2} = \begin{pmatrix}
p x_1\otimes x_1 & q x_2\otimes x_1 + (p-q) x_1\otimes x_2 \\
p x_1\otimes x_2 & -q x_2\otimes x_2 + k x_1\otimes x_1
\end{pmatrix}.
\end{align*}

\begin{rem}\label{rem:GGi-dos}
The braidings of the braided vector spaces $ V_{i,1} $ and $ V_{i,3} $ in \cite[Prop. 3.4]{GGi} belong to this case. If $i=0 $ or $i=2 $, then 
the braidings fit directly. But if $i=1 $ or $i=3 $, then the braidings have the shape
\begin{align*}
\widetilde{\hit}_{1,2}=\begin{pmatrix}
			-q & 0 & 0 & k \\ 
			0 & q & p-q & 0 \\ 
			0 & 0 & p & 0 \\ 
			0 & 0 & 0 & p
			\end{pmatrix}
\end{align*}
which is equivalent to $ \hit_{1,2} $ using the equivalences \ref{item:hiet-a} and \ref{item:hiet-b} in \S \ref{subsec:intro_with_Hietarinta}, namely
\begin{align*}
\hit_{1,2} \stackrel{\ref{item:hiet-a}}{\longleftrightarrow} (\hit_{1,2})^t \stackrel{\ref{item:hiet-b}}{\longleftrightarrow}  \widetilde{\hit}_{1,2}.
\end{align*}
\end{rem}
\begin{table}[ht]
	\caption{Nichols algebras of type $ \hit_{1, 2} $}\label{tab:h12}
	\begin{center}
		\begin{tabular}{| c | c |  c | c | c | c |}\hline
			
			$p$ & $-q$ & $ k $ & $ \J(V)$ & Basis & $ \GK $
			\\\hline 
			$-1$ & $ -1 $ & $ 0 $ &$ \lg x_1^2, x_2^2, r_1 \rg$ & $ (*_{1, 2})_1 $ & $ 0$, $\dim  = 4 $
			\\\cline{2-6} 
			 & $ \in\G_N', N\geq 3 $ &  &$ \lg x_1^2, x_2^N, r_1 \rg$ & $ (*_{1, 2})_1 $ & $ 0$, $\dim = 2N $
			\\\cline{2-6} 
			 & \multicolumn{1}{r}{otherwise}&  &$ \lg x_1^2, r_1 \rg$ & $ (*_{1, 2})_2 $ &  $1$
			\\\hline 
			$\in\G_M', M\geq 3 $& $-1$ & & $ \lg x_1^M, r_1, r_2 \rg$ & $ (*_{1, 2})_3 $ & $ 0$, $\dim  = 2M$
			\\\cline{1-1}\cline{3-6}
			$\notin\G_\infty' $&  & & $ \lg  r_1, r_2 \rg$ & $ (*_{1, 2})_2 $ & $1$
			\\\hline
		\end{tabular}
	\end{center}
\end{table}

\begin{pro}\label{prop:case1,2}
If $ p \neq -1 $ and $ q \neq 1 $, then there are no quadratic relations. Otherwise, the Nichols algebras are as in Table \ref{tab:h12}, where

\begin{align*}
r_1 &= x_1x_2 - qx_2x_1, & r_2 &= x_2^2 - \frac{k}{p+1}x_1^2.
\end{align*}
\begin{align*}
&\hspace{-1.42cm}\text{and} & (*_{1, 2})_1 &=  \{x_1^{a_1}x_2^{a_2}: 0\leq a_1 \leq 1, \, 0\leq a_2 < N = \ord(-q) \} ;\\ 
& & (*_{1, 2})_2 &= \{x_1^{a_1}x_2^{a_2}: 0\leq a_1 \leq 1, \, 0\leq a_2 < \infty\}; \\
& & (*_{1, 2})_3 &= \{x_1^{a_1}x_2^{a_2}: 0\leq a_1 < M = \ord p, \, 0\leq a_2 \leq 1 \}.
\end{align*}
\end{pro}

\pf 
Set $u = \lambda_1x_1^2+\lambda_2x_1x_2+\lambda_3x_2x_1+\lambda_4x_2^2$. Then
\begin{align*}
\Delta(u) &= u\ot 1 + 1\ot u + (\lambda_1(1+p)+\lambda_4k)x_1\ot x_1 + \lambda_4(1-q)x_2\ot x_2 \\
&+ (\lambda_2q+\lambda_3)x_2\ot x_1 + (\lambda_2(1+p-q)+\lambda_3p)x_1\ot x_2.
\end{align*}
Then all the assertions on quadratic relations hold. In particular, row $1$ is established. 
To simplify the exposition, we consider two cases:
\begin{enumerate}[leftmargin=*,label=\rm{(\roman*)}] 
\item \label{case1:12} $ p=-1 $;
\item \label{case2:12} $ q = 1 $.
\end{enumerate}
Case \ref{case1:12}: 
We start by the braided Hopf algebra $\widetilde{\toba} = T(V) / \lg  x_1^2, r_1\rg$, 
that covers $\toba(V)$ and is linearly generated by
\begin{align*}
(*_{1, 2})_2 = \{x_1^{a_1}x_2^{a_2}: 0\leq a_1 \leq 1, \, 0\leq a_2 < \infty \}.
\end{align*}

We compute the derivations of the elements in $(*_{1, 2})_2$; we have for $n\geq 2$
\begin{align*}
\partial_1(x_2^n)&= \ta_nx_2^{n-2}x_1, & \partial_2(x_2^n)&= (n)_{-q}\, x_2^{n-1};
\end{align*}
here the coefficients $\ta_n$, found explicitly with \texttt{WolframAlpha}, are
\begin{align}\label{eq1ofcase1of12}
\ta_n=\begin{cases}
k\dfrac{(-1)^nq^{2n-1}+(q-1)q^{n-1}-(-1)^n}{(q-1)(q+1)^2}, & \textrm{if } q^2\neq 1;\\ \\
k\dfrac{(-1)^n n(n-1)}{2}, & \textrm{if } q= -1;\\ \\
k\dfrac{(-1)^n(2n-1)+1}{4}, & \textrm{if } q=1.
\end{cases}
\end{align}
On the other hand, we have for $n\geq 2$
\begin{align*}
\partial_1(x_2^nx_1)&= (-1)^{n-1} x_2^{n-1}, & \partial_2(x_2^nx_1)&= (n-1)_{-q}\, x_2^{n-2}x_1.
\end{align*}

We claim that $(*_{1, 2})_2$ is linearly independent in $\toba(V)$. Otherwise,
we pick a non-trivial relation  $ r = c_1x_2^n + c_2x_2^{n-1}x_1 $  with $ n\geq 3 $ minimal. Then
\begin{align*}
0 &= \partial_1(r) = c_1\ta_n x_2^{n-2}x_1 +c_2(-1)^{n-1} x_2^{n-1}
\end{align*}
what implies $ c_2 = 0 $. Therefore, we get a new relation $ r=x_2^n $ iff 
\begin{align}\label{eq:h12-minimal}
 \ta_n &= 0 & & \text{ and }  & (n)_{-q} &=0.  
\end{align}
Hence, row $3$ of Table \ref{tab:h12} holds.
Next, suppose that $ -q\in \G'_\infty $.  Observe that
\begin{itemize}
\item If $ q=1 $, then, by \eqref{eq1ofcase1of12}, $ \ta_n = 0 $ only if $ k=0 $,  but this is row $1$.
\item If $N =\ord(-q) \geq 3 $, then it follows that $ \ta_N = 0 $ by \eqref{eq1ofcase1of12}.
\end{itemize}
By \eqref{eq:h12-minimal} and arguing as in previous cases, we conclude that $(*_{1, 2})_1$ is a basis of $\toba(V)$.
Hence, row  $2$ of  Table \ref{tab:h12} holds.

\smallbreak

Case \ref{case2:12}: Again we start by the braided Hopf algebra $\widetilde{\toba} = T(V) / \lg  r_1, r_2\rg$, 
that covers $\toba(V)$ and is linearly generated by
\begin{align*}
(*_{1, 2})_4 = \{x_1^{a_1}x_2^{a_2}: 0\leq a_1 < \infty, \, 0\leq a_2 \leq 1\}.
\end{align*}
 It is easy to see that for $n\geq 2$
\begin{align*}
\partial_1(x_1^n)&= (n)_p \, x_1^{n-1}, &\partial_2(x_1^n)&= 0, \\
\partial_1(x_2x_1^n)&= p(n-1)_p \, x_2x_1^{n-2}, &\partial_2(x_2x_1^n)&= x_1^{n-1}.
\end{align*}
Let $ r = c_1x_1^n + c_2x_2x_1^{n-1} $ be a relation in $\toba(V)$ with $ n\geq 3 $ minimal. Then
\begin{align*}
0 &= \partial_2(r) = c_2 x_1^{n-1}.
\end{align*}
what gives us $ c_2 = 0 $. Hence,
\begin{align*}
0 &= \partial_1(r) = c_1(n)_p \, x_1^{n-1}.
\end{align*}
and we have the relation $ x_1^n $ iff $ p\in \G_n', $  $n\geq 2$; but $ p = -1 $ 
is excluded since it turns out to be case \ref{case1:12}.
The rest of the proof goes as in the case  \ref{case1:12}.
\epf

We next discuss the Nichols algebras arising from the equivalence in \cite{H}. 
Here the equivalences \ref{item:hiet-a}, \ref{item:hiet-c} and \ref{item:hiet-a} $\circ$ \ref{item:hiet-c} give different braidings. 

\medbreak
First, \ref{item:hiet-a} gives rise to the braiding
\begin{align*}
(c'(x_i \ot x_j))_{i, j \in \I_2} = \begin{pmatrix}
p x_1\otimes x_1 + k x_2\otimes x_2 & q x_2\otimes x_1  \\
p x_1\otimes x_2 + (p-q) x_2\otimes x_1 & -q x_2\otimes x_2 
\end{pmatrix}.
\end{align*}

\begin{table}[ht]
	\caption{Nichols algebras of type $ \hit_{1, 2} \ (a) $}\label{tab:h12a}
	\begin{center}
		\begin{tabular}{| c | c | c | c | c | c |}\hline
			$p$ & $-q$ & $ k $ & $ \J(V)$ & Basis & $ \GK $
			\\\hline 
			$-1$ & $-1$ & $0$ &$ \lg x_1^2, x_2^2, r_1 \rg$ & $ (*_{1, 2})_1 $ & $ 0$, $\dim  = 4$
			 			\\ \cline{3-6} 
			 &   & $\neq 0$  &$ \lg r_1, r_2 \rg$ & $ (*_{1, 2})_4 $ &  $1$
			\\\cline{2-6} 
			 & $ \in\G_N', N\geq 3 $ &  &$ \lg  x_2^N, r_1, r_2 \rg$ & $ (*_{1, 2})_1 $ & $ 0$, $\dim  = 2N$
			 	\\ \cline{2-6} 
			 & $ \notin\G'_{\infty}$ &  &$ \lg r_1, r_2 \rg$ & $ (*_{1, 2})_2 $ &  $1$
			\\\hline 
			$\in\G_M', M\geq 3 $& $-1$ & & $ \lg x_1^M, x_2^2, r_1 \rg$ & $ (*_{1, 2})_3 $ & $ 0$, $\dim  = 2M$
			\\\cline{1-1}\cline{3-6}
			$\notin\G_\infty' $&  & & $ \lg x_2^2, r_1 \rg$ & $ (*_{1, 2})_4 $ & $1$
			\\\hline
		\end{tabular}
	\end{center}
\end{table}

\begin{pro}   If $ p \neq -1 $ and $ q \neq 1 $, then there are no quadratic relations. Otherwise, the Nichols algebras are as in Table \ref{tab:h12a}, where
\begin{align*}
r_1 &= x_2x_1 - px_1x_2, & r_2 &= (1-q)x_1^2 - kx_2^2.
\end{align*}
\begin{align*}
&\hspace{-1.42cm}\text{and} & (*_{1, 2})_1 &=  \{x_1^{a_1}x_2^{a_2}: 0\leq a_1 \leq 1, \, 0\leq a_2 < N = \ord(-q) \} ;\\ 
& & (*_{1, 2})_2 &= \{x_1^{a_1}x_2^{a_2}: 0\leq a_1 \leq 1, \, 0\leq a_2 < \infty\}; \\
& & (*_{1, 2})_3 &= \{x_1^{a_1}x_2^{a_2}:  0\leq a_1 < M = \ord p, \, 0\leq a_2 \leq 1 \}; \\
& & (*_{1, 2})_4 &= \{x_1^{a_1}x_2^{a_2}:  0\leq a_1 < \infty, \, 0\leq a_2 \leq 1 \}.
\end{align*}
\end{pro}

\pf Similar to the proof of Proposition \ref{prop:case1,2}.
\epf

Second, \ref{item:hiet-c} gives rise to the braiding
\begin{align*}
(c''(x_i \ot x_j))_{i, j \in \I_2} = \begin{pmatrix}
p x_1\otimes x_1  & p x_2\otimes x_1  \\
q x_1\otimes x_2 + (p-q) x_2\otimes x_1 & -q x_2\otimes x_2+ k x_1\otimes x_1 
\end{pmatrix}.
\end{align*}

\begin{table}[ht]
	\caption{Nichols algebras of type $ \hit_{1, 2} \ (c) $}\label{tab:h12c}
	\begin{center}
		\begin{tabular}{| c | c | c | c | c | c |}\hline
			$p$ & $-q$ & $ k $ & $ \J(V)$ & Basis & $ \GK $
			\\\hline
			$-1$ & $ -1 $ & $ 0 $ &$ \lg x_1^2, x_2^2, r_1 \rg$ & $ (*_{1, 2})_1 $ & $ 0$, $\dim  = 4$
			\\ \cline{2-6}
			 & $ \in\G_N', N\geq 3 $ &  &$ \lg x_1^2, x_2^N, r_1 \rg$ & $ (*_{1, 2})_1 $ & $ 0$, $\dim  = 2N$
			\\\cline{2-6}
			 & \multicolumn{1}{r}{otherwise} &  &$ \lg x_1^2, r_1 \rg$ & $ (*_{1, 2})_2 $ &  $1$
			\\\hline 
			$\in\G_M', M\geq 3 $& $-1$ & & $ \lg x_1^M, r_1, r_2 \rg$ & $ (*_{1, 2})_3 $ & $ 0$, $\dim  = 2M$
			\\\cline{1-1}\cline{3-6}
			$\notin\G_\infty' $&  & & $ \lg r_1, r_2 \rg$ & $ (*_{1, 2})_4 $ & $1$
			\\\hline
		\end{tabular}
	\end{center}
\end{table}

\begin{pro}   If $ p \neq -1 $ and $ q \neq 1 $, then there are no quadratic relations. Otherwise, the Nichols algebras are as in Table \ref{tab:h12c}, where
\begin{align*}
r_1 &= x_2x_1 - qx_1x_2, & r_2 &= x_2^2 - \dfrac{k}{p+1}x_1^2.
\end{align*}
\begin{align*}
&\hspace{-1.42cm}\text{and} & (*_{1, 2})_1 &=  \{x_1^{a_1}x_2^{a_2}: 0\leq a_1 \leq 1, \, 0\leq a_2 < N = \ord(-q) \} ;\\ 
& & (*_{1, 2})_2 &= \{x_1^{a_1}x_2^{a_2}: 0\leq a_1 \leq 1, \, 0\leq a_2 < \infty\}; \\
& & (*_{1, 2})_3 &= \{x_1^{a_1}x_2^{a_2}:  0\leq a_1 < M = \ord p, \, 0\leq a_2 \leq 1 \}; \\
& & (*_{1, 2})_4 &= \{x_1^{a_1}x_2^{a_2}:  0\leq a_1 < \infty, \, 0\leq a_2 \leq 1 \}.
\end{align*}
\end{pro}

\pf Similar to the proof of Proposition \ref{prop:case1,2}.
\epf

Finally, \ref{item:hiet-a} composed with \ref{item:hiet-c} gives rise to the braiding
\begin{align*}
(c'''(x_i \ot x_j))_{i, j \in \I_2} = \begin{pmatrix}
p x_1\otimes x_1 + k x_2\otimes x_2 & p x_2\otimes x_1 + (p-q) x_1\otimes x_2 \\
q x_1\otimes x_2  & -q x_2\otimes x_2 
\end{pmatrix}.
\end{align*}

\begin{table}[ht]
	\caption{Nichols algebras of type $ \hit_{1, 2} \ (ac) $}\label{tab:h12ac}
	\begin{center}
		\begin{tabular}{| c | c | c | c | c | c |}\hline
			$p$ & $-q$ & $ k $ & $ \J(V)$ & Basis & $ \GK $
			\\\hline 
			$-1$ & $-1$ & $0$ &$ \lg x_1^2, x_2^2, r_1 \rg$ & $ (*_{1, 2})_1 $ & $ 0$, $\dim  = 4$
			 			\\ \cline{3-6} 
			 &   & $\neq 0$  &$ \lg r_1, r_2 \rg$ & $ (*_{1, 2})_4 $ &  $1$
			\\\cline{2-6} 
			 & $ \in\G_N', N\geq 3 $ &  &$ \lg  x_2^N, r_1, r_2 \rg$ & $ (*_{1, 2})_1 $ & $ 0$, $\dim  = 2N$
			 	\\ \cline{2-6} 
			 & $ \notin\G'_{\infty}$ &  &$ \lg r_1, r_2 \rg$ & $ (*_{1, 2})_2 $ &  $1$
			\\\hline 
			$\in\G_M', M\geq 3 $& $-1$ & & $ \lg x_1^M, x_2^2, r_1 \rg$ & $ (*_{1, 2})_3 $ & $ 0$, $\dim  = 2M$
			\\\cline{1-1}\cline{3-6}
			$\notin\G_\infty' $&  & & $ \lg x_2^2, r_1 \rg$ & $ (*_{1, 2})_4 $ & $1$
			\\\hline
		\end{tabular}
	\end{center}
\end{table}

\begin{pro}   If $ p \neq -1 $ and $ q \neq 1 $, then there are no quadratic relations. Otherwise, the Nichols algebras are as in Table \ref{tab:h12ac}, where
\begin{align*}
r_1 &= x_1x_2 - px_2x_1, & r_2 &= (1-q)x_1^2 - kx_2^2.
\end{align*}
\begin{align*}
&\hspace{-1.42cm}\text{and} & (*_{1, 2})_1 &=  \{x_1^{a_1}x_2^{a_2}: 0\leq a_1 \leq 1, \, 0\leq a_2 < N = \ord(-q) \} ;\\ 
& & (*_{1, 2})_2 &= \{x_1^{a_1}x_2^{a_2}: 0\leq a_1 \leq 1, \, 0\leq a_2 < \infty\} \\
& & (*_{1, 2})_3 &= \{x_1^{a_1}x_2^{a_2}:  0\leq a_1 < M= \ord p, \, 0\leq a_2 \leq 1 \}; \\
& & (*_{1, 2})_4 &= \{x_1^{a_1}x_2^{a_2}:  0\leq a_1 < \infty, \, 0\leq a_2 \leq 1 \}.
\end{align*}
\end{pro}

\pf Similar to the proof of Proposition \ref{prop:case1,2}.
\epf

\subsection{Case $ \hit_{1, 3}$}\label{subsec:hit13} We assume that $ k\neq 0$, and either  $p\neq 0$, or $q\neq 0$. The associated braiding is $ (c(x_i \ot x_j))_{i, j \in \I_2} = $
\begin{align*}
= \begin{pmatrix}
k^2 x_1\otimes x_1 & k^2 x_2\otimes x_1 -kp x_1\otimes x_1 \\
kx_1\otimes(k x_2 + p x_1) & k^2 x_2\otimes x_2 + pq x_1\otimes x_1 + kq x_1\otimes x_2 -kq x_2\otimes x_1
\end{pmatrix}.
\end{align*}

\begin{table}[ht]
	\caption{Nichols algebras of type $ \hit_{1, 3} $}\label{tab:h13}
	\begin{center}
		\begin{tabular}{|c | c | c | c |}\hline
			$k^2$ & $ \J(V)$ & Basis & $ \GK $
			\\\hline 
			$-1$ & $ \lg x_1^2, x_2^2-kqx_1x_2, x_1x_2 + x_2x_1 \rg$ & $ (*_{1, 3})_1 $ & $ 0$, $\dim  = 4$
			\\\hline
			$1$ & $ \lg x_1x_2 - x_2x_1 +kpx_1^2 \rg$ & $ (*_{1, 3})_2 $ & $2$
			\\\hline
		\end{tabular}
	\end{center}
\end{table}

\begin{pro}   If $ k^4 \neq 1 $, then there are no quadratic relations.
Otherwise, the Nichols algebras are as in Table \ref{tab:h13}, where
\begin{align*}
(*_{1, 3})_1 &=  \{x_1^{a_1}x_2^{a_2}: 0\leq a_i \leq 1 \}, &  (*_{1, 3})_2 &= \{x_1^{a_1}x_2^{a_2}: 0\leq a_i < \infty\}. 
\end{align*}
\end{pro}

\pf  
Set $u = \lambda_1x_1^2+\lambda_2x_1x_2+\lambda_3x_2x_1+\lambda_4x_2^2$. Then
\begin{align*}
\Delta(u) &= u\ot 1 + 1\ot u + (\lambda_1(1+k^2)-kp(\lambda_2-\lambda_3)+\lambda_4pq)x_1\ot x_1 \\
&+ \lambda_4(1+k^2)x_2\ot x_2 + (\lambda_2+\lambda_3k^2+\lambda_4kq)x_1\ot x_2 \\
&+ (\lambda_2k^2+\lambda_3-\lambda_4kq)x_2\ot x_1.
\end{align*}
Thus the quadratic relations above are clear; hence, row $1$ is established. 
By a routine argument, $(*_{1, 3})_2$ generates linearly $\bB := T(V)/ \langle x_1x_2 - x_2x_1 +kpx_1^2 \rangle$. As in the (proof of) Proposition \ref{prop:Nichols23} case \ref{case3:23}, we set 
\begin{align*}
\partial_2(x_2^n) = \sum_{0\le j\le n-1} d_{j}^{(n-1)} x_1^jx_2^{n-1-j}.
\end{align*}
In a similar way as for \eqref{eq:3of2.3}, we prove that 
\begin{align}\label{eq:3of1.3}
d_{ 0}^{(n-1)}=n, \, n\geq 1.
\end{align}
Also observe that, by induction, we have for $ n\geq 1 $
\begin{align}\label{eq:1ofclaim:1.3}
c(x_1^n\ot x_2)= x_2\ot x_1^n -nkp x_1\ot x_1^n.
\end{align}
Then, by \eqref{eq:1ofclaim:1.3}, for $ n\geq 3 $ and $ 0 < i \leq n $,
\begin{align}\label{eq:1of1.3}
\begin{split}
\partial_1(x_1^ix_2^{n-i})&= ix_1^{i-1}x_2^{n-i}+x_1^i\partial_1(x_2^{n-i})-ikpx_1^i\partial_2(x_2^{n-i}),\\
\partial_2(x_1^ix_2^{n-i})&= x_1^{i}\partial_2(x_2^{n-i}).
\end{split}
\end{align}

Assume that the image of $B$ under the projection $\pi: \bB\to \toba(V)$ is not linearly independent. 
Pick $r = \sum_{i=0}^n c_i x_1^ix_2^{n-i}$ a non-trivial homogeneous  relation of minimal degree $ n>2$. By \eqref{eq:1of1.3},
\begin{align}\label{eq:2of1.3}
0 = \partial_2(r) = \sum_{0\le i\le n-1} c_i x_1^{i}\partial_2(x_2^{n-i}).
\end{align}
Looking at the monomials $ x_2^{n-i} $in  \eqref{eq:2of1.3} as in case \ref{case3:23} of the proof of Proposition \ref{prop:Nichols23}, we get $ c_i = 0, \, 0\leq i <n $. Also $ \partial_1(x_1^n) =  nx_1^{n-1} $, then $ r = 0 $. 
\epf

We next discuss the Nichols algebras arising from the equivalence in \cite{H}. 
First, \ref{item:hiet-c} gives rise to the braiding $ (c'(x_i \ot x_j))_{i, j \in \I_2} $
\begin{align*}
 = \begin{pmatrix}
k^2 x_1\otimes x_1 & k^2 x_2\otimes x_1 +kp x_1\otimes x_1 \\
kx_1\otimes(k x_2 - p x_1) & k^2 x_2\otimes x_2 + pq x_1\otimes x_1 - kq x_1\otimes x_2 +kq x_2\otimes x_1
\end{pmatrix}.
\end{align*}
But this is $ \hit_{1, 3}$ up to $p \mapsto -p$, $q \mapsto -q$, so no new Nichols algebra arises.

Second, \ref{item:hiet-a} gives rise to the braiding $ (c''(x_i \ot x_j))_{i, j \in \I_2} $
\begin{align*}
 = \begin{pmatrix}
k^2 x_1\otimes x_1 + pq x_2\otimes x_2 + kp x_1\otimes x_2 -kp x_2\otimes x_1 & k x_2\otimes (kx_1 - qx_2) \\
k^2x_1\otimes x_2 +kq x_2\otimes x_2 & k^2 x_2\otimes x_2 
\end{pmatrix}.
\end{align*}
But this is \ref{item:hiet-b} up to $p \mapsto -q$, $q \mapsto -p$, so no new Nichols algebra arises. Therefore, \ref{item:hiet-a} composed with \ref{item:hiet-c} is \ref{item:hiet-b} up to $p \mapsto q$, $q \mapsto p$, so again no new Nichols algebra arises.

\subsection{Case $ \hit_{1, 4}$}\label{subsec:hit14} We assume that $ k, p, q\neq 0 $. The associated braiding is
\begin{align*}
(c(x_i \ot x_j))_{i, j \in \I_2} = \begin{pmatrix}
q x_2\otimes x_2 & k x_1\otimes x_2 \\
k x_2\otimes x_1 & p x_1\otimes x_1
\end{pmatrix}.
\end{align*}

\begin{table}[ht]
	\caption{Nichols algebras of type $ \hit_{1, 4} $}\label{tab:h14}
	\begin{center}
		\begin{tabular}{| c | c | c |c | c |}\hline
			$k$ & $pq$ &  $ \J(V)$ & Basis & $ \GK $
			\\\hline 
			$-1$ & $ \in\G'_N, N\geq 1 $   &$ \lg  r_{1, 2},r_{2, 2},r_{3, N} \rg$ & $ (*_{1, 4})_1 $ & $ 0$, $\dim  = 4N$
			\\\cline{2-5}
			& $\not\in\G_\infty$   & $ \lg  r_{1, 2},r_{2, 2} \rg$  & $ (*_{1, 4})_2 $ & $1$
			\\\hline
			$\in\G'_N, N\geq 3 $ & $1$  &$ \lg r_{1, N},r_{2, N},r_{3, 1} \rg$ & $ (*_{1, 4})_3 $ & $ 0$, $\dim  = N^2$
			\\\cline{1-1}\cline{3-5}
			$\not\in\G'_\infty$ &  &  $ \lg r_{3, 1} \rg$ &$ (*_{1, 4})_4 $ & $2$
			\\\hline
		\end{tabular}
	\end{center}
\end{table}

\begin{pro} If $ k \neq -1 $ and $ pq \neq 1 $, then there are no quadratic relations. Otherwise, the Nichols algebras are as Table \ref{tab:h14}, where
\begin{align}\label{eqn:relation1of1.4}
r_{1,n}&=\begin{cases}
 (x_2x_1)^{\frac{n}{2}}, & \textrm{if } n \textrm{ is even;}\\
x_1(x_2x_1)^{\frac{n-1}{2}}, & \textrm{if } n \textrm{ is odd.} \\
\end{cases}\\\label{eqn:relation2of1.4}
r_{2,n}&=\begin{cases}
 x_1(x_2x_1)^{\frac{n-2}{2}}x_2, & \textrm{if } n \textrm{ is even;}\\
(x_2x_1)^{\frac{n-1}{2}}x_2, & \textrm{if } n \textrm{ is odd.} \\
\end{cases}\\\label{eqn:relation3of1.4}
r_{3,n}&= x_1^{2n} + (-q)^nx_2^{2n}.
\end{align}
The bases are
\begin{align*}
(*_{1, 4})_1 &=  \{x_i^{a}: 0\leq a < 2(N+1) - i \};\\ 
 (*_{1, 4})_2 &= \{x_i^{a}: 0\leq a < \infty, \, i\in \I_2 \}; \\
 (*_{1, 4})_3 &= \{x_1^{a}(x_2x_1)^{b}x_2^{c}: 0\leq a, b, c; \, c \leq 1, \, S:=a+2b+c\leq 2N-2\\
& \hspace{3cm}\text{ and if } S\geq N  \text{ then } a\geq  S-N+2 \}; \\
 (*_{1, 4})_4 &= \{x_1^{a}(x_2x_1)^{b}x_2^{c}: 0\leq c \leq 1, 0\leq a, b < \infty \}.
\end{align*}
\end{pro}

\pf Set $u = \lambda_1x_1^2+\lambda_2x_1x_2+\lambda_3x_2x_1+\lambda_4x_2^2$. Then
\begin{align*}
\Delta(u) &= u\ot 1 + 1\ot u + (\lambda_1+\lambda_4p)x_1\ot x_1 + (\lambda_4+\lambda_1q) x_2\ot x_2 \\
&+ \lambda_2(1+k)x_1\ot x_2 + \lambda_3(1+k)x_2\ot x_1.
\end{align*}
Then all assertions made about quadratic relations hold. To simplify, we consider two cases:
\begin{enumerate}[leftmargin=*,label=\rm{(\roman*)}] 
\item \label{case1:14} $ k = -1 $;
\item \label{case2:14} $ k\neq -1 $ and $ pq=1 $.
\end{enumerate}
Case \ref{case1:14}:
We consider the braided Hopf algebra 
\begin{align*}
\widetilde{\toba} = T(V) / \lg  r_{1,2}, r_{2,2}\rg = T(V) / \lg  x_1 x_2, x_2 x_1 \rg;
\end{align*}
clearly it is linearly generated by
\begin{align*}
(*_{1, 4})_2 =  \{x_i^{a}: 0\leq a < \infty, \, i\in \I_2 \};
\end{align*}
Note that, for $ n\geq 1 $
\begin{align*}
c(x_1^n\ot x_1)=\begin{cases}
(-q)^{\frac{n}{2}} x_1\ot x_2^n, & \textrm{if } n \textrm{ is even;}\\
-(-q)^{\frac{n+1}{2}}x_2\ot x_2^n, & \textrm{if } n \textrm{ is odd.}
\end{cases}\\
c(x_2^n\ot x_2)=\begin{cases}
(-p)^{\frac{n}{2}} x_2\ot x_1^n, & \textrm{if } n \textrm{ is even;}\\
-(-p)^{\frac{n+1}{2}}x_1\ot x_1^n, & \textrm{if } n \textrm{ is odd.}
\end{cases}
\end{align*}
Then, for $n\geq 2$, we prove that
\begin{align*}
&\begin{split}
\partial_1(x_1^n)&=x_1^{n-1},\\
\partial_2(x_2^n)&=x_2^{n-1}, \\
\partial_1(x_1^n)&=x_1^{n-1}+(-q)^{\frac{n-1}{2}} x_2^{n-1}, \\
\partial_2(x_2^n)&=x_2^{n-1}+(-p)^{\frac{n-1}{2}} x_1^{n-1},
\end{split} 
&& \begin{split} 
\partial_2(x_1^n)&=-(-q)^{\frac{n}{2}} x_2^{n-1},\\
\partial_1(x_2^n)&=-(-p)^{\frac{n}{2}} x_1^{n-1},\\
\partial_2(x_1^n)&=0,\\
\partial_1(x_2^n)&=0,
\end{split}
\hspace{0.1cm}
\begin{split}
&\textrm{\Bigg\} if } n \textrm{ is even;}\\
&\textrm{\Bigg\} if } n \textrm{ is odd.}
\end{split}
\end{align*}
Let $0 \neq r = c_1x_1^n + c_2x_2^{n} \in \ker (\widetilde{\toba} \to \toba(V)) $ with $ n\geq 3 $ minimal. 
If $ n $ is odd, then it is clear that $r = 0$ (just apply $\partial_1$ and $\partial_2$). 
So, suppose that $n$ is even. Then
\begin{align*}
0 &= \partial_1(r) = c_1 x_1^{n-1} -(-p)^{\frac{n}{2}}c_2x_1^{n-1} = (c_1 -(-p)^{\frac{n}{2}}c_2)x_1^{n-1}, \\
0 &= \partial_2(r) = -(-q)^{\frac{n}{2}}c_1 x_2^{n-1}+c_2x_2^{n-1} = (-(-q)^{\frac{n}{2}}c_1 +c_2)x_2^{n-1}.
\end{align*}
This implies that 
\begin{align}\label{eq:h14-system}
0 &= c_1 -(-p)^{\frac{n}{2}}c_2, &0 &= -(-q)^{\frac{n}{2}}c_1 +c_2.
\end{align}
The system above has a non trivial solution iff $ (pq)^{\frac{n}{2}} = 1 $. Hence, if $pq\notin \G_{\infty}$, then $\widetilde{\toba} \simeq \toba(V)$ and arguing as before, we conclude that $(*_{1, 4})_2$ is a basis of $\toba(V)$. This establishes the second row of Table \ref{tab:h14}.
\medbreak
Next, we assume that $pq\in \G_{\infty}$ and set $N = \ord (pq)$. Then, the relation $ r_{3,N} $ \eqref{eqn:relation3of1.4} holds in $\toba(V)$, and we may consider the braided Hopf algebra 
\begin{align*}
\toba = T(V) / \lg  r_{1,2}, r_{2,2}, r_{3,N}\rg  = T(V) / \lg  x_1 x_2, x_2 x_1, x_1^{2N} + (-q)^N x_2^{2N} \rg;
\end{align*}
clearly $(*_{1, 4})_1$ generates linearly $\toba$. 
But $(*_{1, 4})_1$ is linearly independent in $\toba(V)$,
as any linear relation would appear in a degree higher than $\deg r_{3,N}$  
by the definition of $N$ and the discussion above on \eqref{eq:h14-system}. 
This establishes the first row of Table \ref{tab:h14}.

\smallbreak
Case \ref{case2:14}: We consider the braided Hopf algebra 
\begin{align*}
\widetilde{\toba} = T(V) / \lg  r_{3,1}\rg = T(V) / \lg  x_1^2 -q x_2^2 \rg;
\end{align*}
We claim that $(*_{1, 4})_4$ linearly generates $\widetilde{\toba}$.

\emph{Proof of the claim. }
It is enough to show that $I = $ linear span of $(*_{1, 4})_4$ is left ideal, since clearly $1 \in I$.
Obviously, $x_1I \subset I$, so we only need to show that $x_2I \subset I$. Observe first that
\begin{align*}
x_2 x_1^2 &= q x_2^3 = x_1^2 x_2, &  x_2(x_2x_1) &= p x_1^3.
\end{align*}
Thus we multiply $x_2 x_1^{a}(x_2x_1)^{b}x_2^{c}$ and get
\begin{align*}
&x_2 x_1^{2h + 1}(x_2x_1)^{b}x_2^{c} =  x_1^{2h}(x_2x_1)^{b + 1}x_2^{c}, &  a&= 2h + 1 \text{ odd;}
\\
&x_2 x_1^{a}x_2^{c} =  x_1^{a} x_2^{c+1} =  p^c x_1^{2(h + c)}x_2^{1-c}, &  a&= 2h \text{ even;}
\\
&x_2 x_1^{a}(x_2x_1)^{b}x_2^{c} = x_1^{a} x_2(x_2x_1)^{b}x_2^{c} = px_1^{a + 3} (x_2x_1)^{b-1}x_2^{c}, &  a&= 2h \text{ even, }
b > 0.
\end{align*}

The claim is proved. Then we check, by an inductive argument, that 
\begin{align*}
\partial_1(x_1^a) &= \begin{cases}
(\frac{a-2}{2})_k \, x_1^{a-1}, \qquad \,\, & \textrm{if } a \textrm{ is even;}\\
(\frac{a-1}{2})_k \, x_1^{a-1}, & \textrm{if } a \textrm{ is odd.}
\end{cases} \\
\partial_2(x_1^a) &= \begin{cases}
q(\frac{a-2}{2})_k \, x_1^{a-2}x_2, & \textrm{if } a \textrm{ is even;}\\
q(\frac{a-3}{2})_k \, x_1^{a-3}x_2x_1, & \textrm{if } a \textrm{ is odd, } a\geq 3;\\
0, & \textrm{if } a =1,
\end{cases}
\end{align*}
$ a\in \N$.
Also, we have that, for $ b\geq 1$
\begin{align*}
\partial_1((x_2x_1)^b) &= 0, && \partial_2((x_2x_1)^b) = (2b-1)_k \, x_1(x_2x_1)^{b-1}.
\end{align*}
Therefore, for $ a,b\geq 1 $ 
\begin{align*}
\partial_1(x_1^a(x_2x_1)^b) &= \begin{cases}
(\frac{a-2}{2})_k \, x_1^{a-1}(x_2x_1)^b, & \textrm{if } a \textrm{ is even;}\\
(\frac{a-1}{2}+2b)_k \, x_1^{a-1}(x_2x_1)^b, & \textrm{if } a \textrm{ is odd.}
\end{cases} \\
\partial_2(x_1^a(x_2x_1)^b) &= \begin{cases}
(\frac{a-2}{2}+2b)_k \, x_1^{a+1}(x_2x_1)^{b-1}, & \textrm{if } a \textrm{ is even;}\\
q(\frac{a-3}{2})_k \, x_1^{a-3}(x_2x_1)^{b+1}, & \textrm{if } a \textrm{ is odd, } a\geq 3;\\
0, & \textrm{if } a =1.
\end{cases}
\end{align*}
and then
\begin{align*}
\partial_1(x_1^a(x_2x_1)^bx_2) &= \begin{cases}
\partial_1(x_1^a(x_2x_1)^b)x_2, & \textrm{if } a \textrm{ is even;}\\
(\frac{a+1}{2}+2b)_k \, x_1^{a-1}(x_2x_1)^bx_2, & \textrm{if } a \textrm{ is odd.}
\end{cases} \\
\partial_2(x_1^a(x_2x_1)^bx_2) &= \begin{cases}
(\frac{a}{2}+2b)_k \, x_1^{a+1}(x_2x_1)^{b-1}x_2, & \textrm{if } a \textrm{ is even;}\\
\partial_2(x_1^a(x_2x_1)^b)x_2, & \textrm{if } a \textrm{ is odd, } a\geq 3;\\
0, & \textrm{if } a =1.
\end{cases}
\end{align*}
We then proceed as in the previous case. Namely, let $r \in \ker (\widetilde{\toba} \to \toba(V))$
be an homogeneous relation of degree $n$ with $n\geq 3$ minimal.
We consider separately the cases $n$ odd and $n$ even.
Using the  derivations, we see that the relations $r_{1,N}$ and $r_{2,N}$ hold.
Also, any other relation arises in higher degree.
In this way, rows $3$ and $4$ of Table  \ref{tab:h14}   are established.
\epf

We finally discuss the Nichols algebras arising from the equivalence in \cite{H}. 
First, \ref{item:hiet-a} and \ref{item:hiet-a} composed with \ref{item:hiet-c} give rise to the same braiding
\begin{align*}
(c'(x_i \ot x_j))_{i, j \in \I_2} = \begin{pmatrix}
q x_2\otimes x_2 & k x_1\otimes x_2 \\
k x_2\otimes x_1 & p x_1\otimes x_1
\end{pmatrix}
\end{align*}
which is $ \hit_{1, 4}$ up to $p \leftrightarrow q$, so no new Nichols algebra arises.

Second, \ref{item:hiet-c} gives the initial $ \hit_{1, 4}$, so no new Nichols algebra arises.

\subsection{Case $ \hit_{0, 1}$}\label{subsec:hit01}  The associated braiding is
\begin{align*}
(c(x_i \ot x_j))_{i, j \in \I_2} = \begin{pmatrix}
kx_1\otimes x_1 & -kx_2\otimes x_1 \\
-kx_1\otimes x_2 & kx_2\otimes x_2 + kx_1\otimes x_1
\end{pmatrix}.
\end{align*}

\begin{pro} 
If $k^2 \neq 1$, then there are no quadratic relations. 

\smallbreak
If $k=1$, then 
\begin{align}\label{eqn:Nichols10,k1}
\toba(V) =T(V)/  \lg  x_1x_2 + x_2x_1 \rg,
\end{align}
$B_1 = \{x_1^{a_1}x_2^{a_2}:  a_i  \in \N_0 \}$  is a PBW-basis of $\toba(V)$ and 
$\GK \toba(V) = 2$.

\smallbreak
If $k=-1$, then 
\begin{align}\label{eqn:Nichols10,k-1}
\toba(V) =T(V)/   \lg  x_1^2, x_1x_2 - x_2x_1 \rg,
\end{align}
$B_2 =  \{x_1^{a_1}x_2^{a_2}: a_1 \in \I_{0,2}, \, a_2 \in \N_0 \}$  is a PBW-basis of $ \toba(V) $. Hence \newline
$\GK \toba(V) = 1$.
\end{pro}

\pf Set $u = \lambda_1x_1^2+\lambda_2x_1x_2+\lambda_3x_2x_1+\lambda_4x_2^2$. Then
\begin{align*}
\Delta(u) &= u\ot 1 + 1\ot u + ((1+k)\lambda_1+k\lambda_4)x_1\ot x_1 +(1+k)\lambda_4 x_2\ot x_2 \\
&+ (\lambda_2-k\lambda_3)x_1\ot x_2 + (\lambda_3-k\lambda_2)x_2\ot x_1.
\end{align*}
From here, all assertions about quadratic relations hold. 

Case $k =1$: Let  $\bB = T(V)/ \langle x_1x_2 + x_2x_1 \rangle$. 
By a standard argument, $B_1$ is a basis of $\bB$. Note that, for $n\geq 2$
\begin{align*}
\partial_1(x_2^n)&= \binom{n}{2} x_1x_2^{n-2}, & \partial_2(x_2^n)&=n x_2^{n-1}.
\end{align*}
Let $ n\geq 3 $ and $i \in \I_{2, n-2}$. Then
\begin{align}\label{eqn:prop:0.1}
\begin{split}
\partial_1(x_1^ix_2^{n-i})&= \binom{n-i}{2} x_1^{i+1}x_2^{n-i-2}+ix_1^{i-1}x_2^{n-i}, \\
\partial_2(x_1^ix_2^{n-i})&=(-1)^i (n-i) x_1^ix_2^{n-i-1}.
\end{split}
\end{align}
Observe that the formulae \eqref{eqn:prop:0.1} also hold for $i\in \{0, 1, n-1, n\} $ taking by $0$ 
the terms that are not well defined. 

Assume that the image of $B_1$ under the projection $\pi: \bB\to \toba(V)$ is not linearly independent. 
Pick $0 \neq r = \sum_{i=0}^N c_i x_1^ix_2^{N-i}$ a homogeneous  relation of minimal degree $ N>2$. Thus
\begin{align*}
0 &= \partial_2(r) = \sum_{i=0}^{N-1} c_i (-1)^i(N-i) x_1^ix_2^{N-i-1}.
\end{align*}
Then, $ c_i = 0 $ for $i \in \I_{0, N-1}$. But $ \partial_1(x_1^N) = N x_1^{N-1} $, hence $r = 0$.

\smallbreak
Case $k = -1$: Let  $\bB = T(V)/ \langle x_1^2, x_1x_2 - x_2x_1 \rangle$. 
By a standard argument, $B_2$ is a basis of $\bB$. Note that, for $n\geq 2$
\begin{align*}
\partial_1(x_2^n) &= \begin{cases}
-\frac{n}{2} x_1x_2^{n-2}, & \textrm{if } n \textrm{ is even;}\\
-\frac{n-1}{2} x_1x_2^{n-2}, & \textrm{if } n \textrm{ is odd.}
\end{cases}, & \partial_2(x_2^n)&=\begin{cases}
0, & \textrm{if } n \textrm{ is even;}\\
x_2^{n-1}, & \textrm{if } n \textrm{ is odd.}
\end{cases}.
\end{align*}
Then, for $n\geq 2$,
\begin{align*}
\partial_1(x_1x_2^n) &= x_2^{n}, & \partial_2(x_1x_2^n)&=\begin{cases}
0, & \textrm{if } n \textrm{ is even;}\\
x_1x_2^{n-1}, & \textrm{if } n \textrm{ is odd.}
\end{cases}.
\end{align*}

Suppose that the image of $B$ under the projection $\pi: \bB\to \toba(V)$ is not linearly independent. 
Pick $r = c_1 x_2^{N} + c_2 x_1x_2^{N-1}$ a homogeneous non-trivial relation of minimal degree $ N>2$. Applying $ \partial_1 $ to $ r $, we obtain $r = 0$.
\epf

We next discuss the Nichols algebras arising from the equivalence in \cite{H}. 
First, \ref{item:hiet-a}, \ref{item:hiet-b} and \ref{item:hiet-a} composed with \ref{item:hiet-c} give rise to the same braiding
\begin{align*}
(c'(x_i \ot x_j))_{i, j \in \I_2} = \begin{pmatrix}
kx_1\otimes x_1 + kx_2\otimes x_2 & -kx_2\otimes x_1 \\
-kx_1\otimes x_2 & kx_2\otimes x_2 
\end{pmatrix}.
\end{align*}
Since  \ref{item:hiet-b} is a change of basis, the Nichols algebras are isomorphic.

Finally, \ref{item:hiet-c} gives the initial braiding
$ \hit_{0, 1}$, so no new Nichols algebra arises.

\section{Appendix}

Here we collect all isomorphism classes of algebras arising as Nichols algebras in Theorem \ref{th:main}, 
see the information in Table \ref{tab:general}. In many cases a change of variables is needed, and we leave to the reader its explicit calculation. All the algebras are of the form $T(W)/ \J$, where $W$ has a basis $y_1, y_2$. 

\newcommand{\letra}{\zeta}
\newcommand{\letrb}{\eta}
\begin{table}[ht]
	\caption{Algebras arising as Nichols algebras of rank 2 }\label{tab:appendix}
	\begin{center}
		\begin{tabular}{|p{3cm} | p{5cm} | p{3cm} | } 
	\hline	Algebra & $\J$ 	& Parameters	\\	\hline
\begin{small}quantum plane \end{small}  &	$\lg y_1y_2 - \letra y_2 y_1 \rg$		&  $\letra \in \ku^{\times}$ 
			\\ \hline
\begin{small}quantum plane \end{small} &	$\lg y_1^M, y_1y_2 - \letra y_2 y_1 \rg$		&  $\letra \in \ku^{\times}$,\newline $M \in \N_{\ge2}$ 
\\ \hline
\begin{small}quantum plane \end{small} &	$\lg y_1^M,  y_1y_2 - \letra y_2 y_1, y_2^N \rg$		&  $\letra \in \ku^{\times}$, \newline $M, N \in \N_{\ge2}$ 
\\ \hline
\begin{small}deformation of\newline a quantum plane \end{small} &	$\lg y_2^2 - \letra y_1^2, y_1y_2 -\letrb y_2 y_1, y_1^N  \rg$		&  $\letrb, \letra \in \ku^{\times}$, \newline
$N \in \N_{\ge 3}$ 
\\ \hline
\begin{small}deformation of\newline a quantum plane \end{small} &	$\lg y_2^2 - \letra y_1^2, y_1y_2 -\letrb y_2 y_1 \rg$		&  $\letrb, \letra \in \ku^{\times}$ 
\\ \hline
\begin{small}deformation of\newline an exterior algebra \end{small}&	$\lg y_1^2,  y_2^2 - y_1 y_2, y_1 y_2 + y_2 y_1 \rg$		&  
\\ \hline
&	$\lg y_1^2 - \letra y_2^2, y_1y_2 + \epsilon y_2 y_1, y_1y_2^N  \rg$		&  $\letra \in \ku^{\times}$, 
 $\epsilon \in \G_2$, \newline $N \in \N_{\ge2}$
\\ \hline
\begin{small}Jordan plane \end{small}  &	$\lg y_1^2 -y_1y_2 + y_2 y_1 \rg$		&   
\\ \hline
\begin{small} super Jordan  plane \end{small}  &	$\lg y_1^2,  y_2^2y_1 -  y_1 y_2 y_1 - y_1 y_2^2 \rg$ &  
\\ \hline
&	$\lg y_1 y_2, y_2 y_1 \rg$		&
\\ \hline
 &	$\lg y_1 y_2, y_2 y_1, y_1^{2N}-\letra y_2^{2N} \rg$		&  $\letra \in \ku^{\times}$,   $N \in \N$ 
\\ \hline
 &	$\lg y_1^2 - \letra y_2^2 \rg$	& $\letra \in \ku^{\times}$ 
\\ \hline
 &	$\lg y_1^2 - \letra y_2^2, r_{1,N}, r_{2,N} \rg$, \newline
 cf. \eqref{eqn:relation1of1.4}, \eqref{eqn:relation2of1.4}		&  $\letra \in \ku^{\times}$,   $N \in \N_{\ge 3}$ 
\\ \hline
		\end{tabular}
	\end{center}
\end{table}

\end{document}